\documentclass{amsart}%
\usepackage{amssymb}
\usepackage{amsfonts}
\usepackage{geometry}
\usepackage{graphicx}
\usepackage{amsmath}%
\providecommand{\U}[1]{\protect\rule{.1in}{.1in}}
\newtheorem{theorem}{Theorem}
\theoremstyle{plain}

\newtheorem{corollary}[theorem]{Corollary}

\newtheorem{definition}[theorem]{Definition}

\newtheorem{lemma}[theorem]{Lemma}
\newtheorem{notation}[theorem]{Notation}

\newtheorem{remark}[theorem]{Remark}

\numberwithin{equation}{section}
\begin{document}
\title[Trilinear Fourier characterization]{A single scale smooth Alpert trilinear characterization of the Fourier
extension conjecture on $\mathbb{P}^{2}$}
\author{Cristian Rios}
\address{University of Calgary\\
Calgary, Alberta, Canada}
\email{crios@ucalgary.ca}
\author{Eric Sawyer}
\address{McMaster University\\
Hamilton, Ontario, Canada}
\email{sawyer@mcmaster.ca}
\thanks{Eric Sawyer's research supported in part by a grant from the National Sciences
and Engineering Research Council of Canada}
\maketitle

\begin{abstract}
We show that for every $0<\delta<1$, the Fourier extension conjecture on the
paraboloid in three dimensions is equivalent to a local single scale smooth
Alpert trilinear inequality of the form,%
\[
\left\Vert \mathcal{E}\mathsf{Q}_{s,U_{1}}^{\eta}f_{1}\ \mathcal{E}%
\mathsf{Q}_{s,U_{2}}^{\eta}f_{2}\ \mathcal{E}\mathsf{Q}_{s,U_{3}}^{\eta}%
f_{3}\right\Vert _{L^{\frac{q}{3}}\left(  B_{\mathbb{R}^{3}}\left(
0,2^{\frac{s}{1-\delta}}\right)  \right)  }\leq C_{\delta,\varepsilon,\nu
}2^{\varepsilon s}\left\Vert f_{1}\right\Vert _{L^{\infty}\left(
U_{1}\right)  }\left\Vert f_{2}\right\Vert _{L^{\infty}\left(  U_{2}\right)
}\left\Vert f_{3}\right\Vert _{L^{\infty}\left(  U_{3}\right)  },
\]
which is a variant of the analogous multiscale trilinear inequality in
\cite{RiSa}, where the smooth Alpert projections $\mathsf{Q}_{s,U_{k}}^{\eta}$
were replaced more generally with $\mathsf{Q}_{s_{k},U_{k}}^{\eta}$ where
$0\leq s_{1}\leq s_{2}\leq s_{3}$ and $s_{2},s_{3}$ (but not necessarily
$s_{1}$) were close to $s$.

\end{abstract}
\tableofcontents

\section{Introduction}

In this paper we prove a variant of the main result from \cite{RiSa}. Let
$\mathcal{E}f\equiv\widehat{\Phi_{\ast}f}$ be the Fourier extension operator
given as the Fourier transform of the pushforward of the measure $f\left(
x\right)  dx$ in $\mathbb{R}^{2}$ to the paraboloid $\mathbb{P}^{2}$ by the
map $\Phi\left(  x_{1},x_{2}\right)  =\left(  x_{1},x_{2},x_{1}^{2}+x_{2}%
^{2}\right)  $. Define
\[
\mathsf{Q}_{s,U}f=\sum_{I\in\mathcal{G}_{s}\left[  U\right]  }\bigtriangleup
_{I;\kappa}^{\eta}f=\sum_{I\in\mathcal{G}_{s}\left[  U\right]  }\left\langle
f,h_{I;\kappa}\right\rangle h_{I;\kappa}^{\eta}%
\]
to be the smooth Alpert pseudoprojection of $f$ at level $s$ in the grid
$\mathcal{G}\left[  U\right]  $. Here $h_{I;\kappa}^{\eta}$ is a smooth Alpert
wavelet with $\kappa$ vanishing moments as introduced in \cite{Saw7}, and also
described below. The Fourier extension conjecture for the paraboloid
$\mathbb{P}^{2}$ in $\mathbb{R}^{3}$ is,%
\begin{equation}
\left\Vert \mathcal{E}f\right\Vert _{L^{q}\left(  \mathbb{R}^{3}\right)  }\leq
C\left\Vert f\right\Vert _{L^{q}\left(  U\right)  },\ \ \ \ \ \text{for all
}q>3\text{ and }f\in L^{q}\left(  U\right)  ,\label{FEC}%
\end{equation}
and we will often refer to this inequality as the \emph{linear} Fourier
extension conjecture to distinguish it from the trilinear variants below.

We say that a triple $\left(  U_{1},U_{2},U_{3}\right)  $ of squares in
$U^{3}$ (where $U\subset B_{\mathbb{R}^{2}}\left(  0,\frac{1}{2}\right)  $) is
$\nu$\emph{-disjoint} if
\begin{equation}
\operatorname*{diam}\left(  U_{i}\right)  \approx\operatorname*{diam}\left(
U_{j}\right)  \approx\nu\text{ and }\operatorname*{dist}\left(  U_{i}%
,U_{j}\right)  \geq\nu,\ \ \ \ \ \text{for all }1\leq i\neq j\leq
3.\label{weak sep}%
\end{equation}
This definition essentially coincides with the definition of weak
transversality in $\mathbb{R}^{3}$ introduced by Muscalu and Oliveira in
\cite{MuOl}, where they establish results on Fourier extension involving
tensor functions.

\begin{theorem}
\label{single main Alpert}Suppose $0<\delta<1$ and $\kappa\in\mathbb{N}$ with
$\kappa>\frac{20}{\delta}$. Then the linear Fourier extension conjecture
(\ref{FEC}) for the paraboloid $\mathbb{P}^{2}$ in $\mathbb{R}^{3} $holds
\emph{if and only if} for every $q>3$ there is $\nu>0$ depending on $q$, such
that the single scale disjoint\emph{\ }trilinear inequality%
\[
\left\Vert \mathcal{E}\mathsf{Q}_{s,U_{1}}^{\eta}f_{1}\ \mathcal{E}%
\mathsf{Q}_{s,U_{2}}^{\eta}f_{2}\ \mathcal{E}\mathsf{Q}_{s,U_{3}}^{\eta}%
f_{3}\right\Vert _{L^{\frac{q}{3}}\left(  B_{\mathbb{R}^{3}}\left(
0,2^{\frac{s}{1-\delta}}\right)  \right)  }\leq C_{\delta,\varepsilon,\nu
}2^{\varepsilon s}\left\Vert f_{1}\right\Vert _{L^{\infty}\left(
U_{1}\right)  }\left\Vert f_{2}\right\Vert _{L^{\infty}\left(  U_{2}\right)
}\left\Vert f_{3}\right\Vert _{L^{\infty}\left(  U_{3}\right)  }\ ,
\]
holds for all $s\in\mathbb{N}$, $\varepsilon>0$, all single scale smooth
Alpert pseudoprojections $\mathsf{Q}_{s,U_{k}}f_{k}=\sum_{I\in\mathcal{G}%
_{s}\left[  U_{k}\right]  }\bigtriangleup_{I;\kappa}^{\eta}f$ with $\kappa$
vanishing moments, and all $\nu$-disjoint triples $\left(  U_{1},U_{2}%
,U_{3}\right)  \in U^{3}$.
\end{theorem}

This theorem says that the linear Fourier extension conjecture is equivalent
to a corresponding disjoint trilinear smooth Alpert conjecture at single
scales, that permits the mild growth constant $2^{\varepsilon s}$ and at the
expense of taking integration over the local ball $B_{\mathbb{R}^{3}}\left(
0,2^{\frac{s}{1-\delta}}\right)  $ (larger than $B_{\mathbb{R}^{3}}\left(
0,2^{s}\right)  $), where $s$ is the scale of the smooth Alpert
pseudoprojection. We will give precise definitions below.

The main tool used in the proof is a modification (as in \cite[proof of
Theorem 3]{RiSa}) of the simplest part of the pigeonholing argument of
Bourgain and Guth \cite{BoGu} that derives linear inequalities from
multilinear inequalities. Two key points should be noted. The first is that
using disjoint instead of transverse coordinate patches, essentially
eliminates the problematic Case 3 in the Bourgain and Guth argument, and
permits a characterization. The second key point is that the parabolic
rescaling of Tao, Vargas and Vega that is used in the proof, applies whenever
the family of functions being tested over is invariant under dyadic parabolic
rescalings - such as (1) the family of bounded functions, and (2) the family
of pseudoprojections of smooth Alpert wavelets at a single scale.

Theorem \ref{single main Alpert} will follow from Theorems \ref{without} and
\ref{main Alpert}\ below.

\subsection{Preliminaries}

Parameterize a patch of the paraboloid $\mathbb{P}^{2}$ in the usual way, i.e.
$\Phi:U\rightarrow\mathbb{P}^{2}$ by
\[
z=\Phi\left(  x\right)  \equiv\left(  x,\left\vert x\right\vert ^{2}\right)
=\left(  x_{1},x_{2},x_{1}^{2}+x_{2}^{2}\right)  ,\ \ \ \ \ x\in U\subset
B_{\mathbb{R}^{2}}\left(  0,\frac{1}{2}\right)  .
\]
For $f\in L^{p}\left(  U\right)  $, define the Fourier extension operator
$\mathcal{E}$\ for the paraboloid in three dimensions by
\[
\mathcal{E}f\left(  \xi\right)  =\mathcal{E}_{U}f\left(  \xi\right)
\equiv\mathcal{F}\left(  \Phi_{\ast}\left[  f\left(  x\right)  dx\right]
\right)  \left(  \xi\right)  =\left[  \Phi_{\ast}\left(  f\left(  x\right)
dx\right)  \right]  ^{\wedge}\left(  \xi\right)  =\int_{U}e^{-i\Phi\left(
x\right)  \cdot\xi}f\left(  x\right)  dx,\ \ \ \ \ \text{for }\xi\in
\mathbb{R}^{3},
\]
where $\mathcal{F}$ is the Fourier transform in $\mathbb{R}^{3}$, and where
$\Phi_{\ast}f=\Phi_{\ast}\left(  f\left(  x\right)  dx\right)  $ denotes the
pushforward of the measure $f\left(  x\right)  dx$ to the paraboloid
$\mathbb{P}^{2}$ under the map $\Phi:U\rightarrow\mathbb{P}^{2}$. Note that
the Fourier restriction and extension conjectures are equivalent to
boundedness of $\mathcal{E}_{U}$ from $L^{q}\left(  U\right)  $ to
$L^{q}\left(  \mathbb{R}^{3}\right)  $ for $q>3$ (except for some boundary
points on the Knapp line).

Now we recall the frame of smooth Alpert wavelets constructed in \cite{Saw7}.

\subsection{Smooth Alpert frames}

Here is the construction from \cite{Saw7} of smooth Alpert projections
$\left\{  \bigtriangleup_{Q;\kappa}^{\eta}\right\}  _{Q\in\mathcal{D}}$ and
corresponding wavelets $\left\{  h_{Q;\kappa}^{a,\eta}\right\}  _{Q\in
\mathcal{D},\ a\in\Gamma_{n}}$ of order $\kappa$ in $n$-dimensional space
$\mathbb{R}^{n}$, provided $\eta>0$ is sufficiently small and $\kappa
\in\mathbb{N}$ is sufficiently large. First, we recall that the (unsmoothed)
Alpert wavelets $\left\{  h_{Q;\kappa}^{a}\right\}  _{a\in\Gamma}$ constructed
in \cite{RaSaWi} are an orthonormal basis for the finite dimensional vector
subspace of $L^{2}$ that consists of linear combinations of the indicators
of\ the children $\mathfrak{C}\left(  Q\right)  $ of $Q$ multiplied by
polynomials of degree at most $\kappa-1$, and such that the linear
combinations have vanishing moments on the cube $Q$ up to order $\kappa-1$:%
\[
L_{Q;k}^{2}\left(  \mu\right)  \equiv\left\{
f=\mathop{\displaystyle \sum }\limits_{Q^{\prime}\in\mathfrak{C}\left(
Q\right)  }\mathbf{1}_{Q^{\prime}}p_{Q^{\prime};k}\left(  x\right)  :\int
_{Q}f\left(  x\right)  x_{i}^{\ell}d\mu\left(  x\right)  =0,\ \ \ \text{for
}0\leq\ell\leq k-1\text{ and }1\leq i\leq n\right\}  ,
\]
where $p_{Q^{\prime};k}\left(  x\right)  =\sum_{\alpha\in\mathbb{Z}_{+}%
^{n}:\left\vert \alpha\right\vert \leq k-1\ }a_{Q^{\prime};\alpha}x^{\alpha}$
is a polynomial in $\mathbb{R}^{n}$ of degree $\left\vert \alpha\right\vert
=\alpha_{1}+...+\alpha_{n}$ at most $\kappa-1$, and $x^{\alpha}=x_{1}%
^{\alpha_{1}}x_{2}^{\alpha_{2}}...x_{n-1}^{\alpha_{n-1}}$. Let $d_{Q;\kappa
}\equiv\dim L_{Q;\kappa}^{2}\left(  \mu\right)  $ be the dimension of the
finite dimensional linear space $L_{Q;\kappa}^{2}\left(  \mu\right)  $.
Moreover, for each $a\in\Gamma_{n}$, we may assume the wavelet $h_{Q;\kappa
}^{a}$ is a translation and dilation of the unit wavelet $h_{Q_{0};\kappa}%
^{a}$, where $Q_{0}=\left[  0,1\right)  ^{n}$ is the unit cube in
$\mathbb{R}^{n}$. The Alpert projection $\bigtriangleup_{Q;\kappa}$ onto
$L_{Q;k}^{2}\left(  \mu\right)  $ is given by $\bigtriangleup_{Q;\kappa}%
f=\sum_{a\in\Gamma_{n}}\left\langle f,h_{Q;\kappa}^{a}\right\rangle
h_{Q;\kappa}^{a}$.

Given a small positive constant $\eta>0$, define a smooth approximate identity
by $\phi_{\eta}\left(  x\right)  \equiv\eta^{-n}\phi\left(  \frac{x}{\eta
}\right)  $ where $\phi\in C_{c}^{\infty}\left(  B_{\mathbb{R}^{n}}\left(
0,1\right)  \right)  $ has unit integral, $\int_{\mathbb{R}^{n}}\phi\left(
x\right)  dx=1$, and vanishing moments of \emph{positive} order less than
$\kappa$, i.e.
\begin{equation}
\int\phi\left(  x\right)  x^{\gamma}dx=\delta_{\left\vert \gamma\right\vert
}^{0}=\left\{
\begin{array}
[c]{ccc}%
1 & \text{ if } & \left\vert \gamma\right\vert =0\\
0 & \text{ if } & 0<\left\vert \gamma\right\vert <\kappa
\end{array}
\right.  .\label{van pos}%
\end{equation}
The \emph{smooth} Alpert `wavelets' were then defined in \cite{Saw7} by%
\[
h_{Q;\kappa}^{a,\eta}\equiv h_{Q;\kappa}^{a}\ast\phi_{\eta\ell\left(
Q\right)  },
\]
and we have for $0\leq\left\vert \beta\right\vert <\kappa$,%
\begin{align*}
&  \int h_{Q;\kappa}^{a,\eta}\left(  x\right)  x^{\beta}dx=\int\phi_{\eta
\ell\left(  I\right)  }\ast h_{Q;\kappa}^{a}\left(  x\right)  x^{\beta}%
dx=\int\int\phi_{\eta\ell\left(  I\right)  }\left(  y\right)  h_{Q;\kappa}%
^{a}\left(  x-y\right)  x^{\beta}dx\\
&  =\int\phi_{\eta\ell\left(  I\right)  }\left(  y\right)  \left\{  \int
h_{Q;\kappa}^{a}\left(  x-y\right)  x^{\beta}dx\right\}  dy=\int\phi_{\eta
\ell\left(  I\right)  }\left(  y\right)  \left\{  \int h_{Q;\kappa}^{a}\left(
x\right)  \left(  x+y\right)  ^{\beta}dx\right\}  dy\\
&  =\int\phi_{\eta\ell\left(  I\right)  }\left(  y\right)  \left\{  0\right\}
dy=0,
\end{align*}
by translation invariance of Lebesgue measure.

There is a linear map $S_{\eta}^{\mathcal{D}}=S_{\kappa,\eta}^{\mathcal{D}}$,
bounded and invertible on all $L^{p}\left(  \mathbb{R}^{2}\right)  $ spaces,
$1<p<\infty$,$\,$such that if we define%
\[
\bigtriangleup_{I;\kappa}^{\eta}f\equiv\left(  \bigtriangleup_{I;\kappa
}f\right)  \ast\phi_{\eta\ell\left(  I\right)  },
\]
then%
\[
\bigtriangleup_{I;\kappa}^{\eta}f\equiv\sum_{a\in\Gamma_{n}}\left\langle
\left(  S_{\eta}^{\mathcal{D}}\right)  ^{-1}f,h_{I;\kappa}^{a}\right\rangle
h_{I;\kappa}^{a,\eta}=\sum_{a\in\Gamma_{n}}\left\langle \left(  S_{\eta
}^{\mathcal{D}}\right)  ^{-1}f,h_{I;\kappa}^{a}\right\rangle S_{\eta
}^{\mathcal{D}}h_{I;\kappa}^{a}=\sum_{a\in\Gamma_{n}}\left(  S_{\eta
}^{\mathcal{D}}\bigtriangleup_{I;\kappa}\left(  S_{\eta}^{\mathcal{D}}\right)
^{-1}\right)  f=\sum_{a\in\Gamma_{n}}\bigtriangleup_{I;\kappa}^{\spadesuit
}f\ ,
\]
where $A^{\spadesuit}$ denotes the commutator $S_{\eta}^{\mathcal{D}}A\left(
S_{\eta}^{\mathcal{D}}\right)  ^{-1}$ of an operator $A$ with $S_{\eta
}^{\mathcal{D}}$.

\begin{theorem}
[\cite{Saw7}]\label{reproducing}Let $n\geq2$ and $\kappa\in\mathbb{N}$ with
$\kappa>\frac{n}{2}$. Then there is $\eta_{0}>0$ depending on $n$ and $\kappa
$\ such that for all $0<\eta<\eta_{0}$, and for all grids $\mathcal{D}$ in
$\mathbb{R}^{n} $, and all $1<p<\infty$, there is a bounded invertible
operator $S_{\kappa,\eta}^{\mathcal{D}}$ on $L^{p}$, and a positive constant
$C_{p,n,\eta}$ such that the collection of functions $\left\{  h_{I;\kappa
}^{a,\eta}\right\}  _{I\in\mathcal{D},\ a\in\Gamma_{n}}$ is a $C_{p,n,\eta}%
$-frame for $L^{p}$, by which we mean,%
\begin{align}
f\left(  x\right)   &  =\sum_{I\in\mathcal{D}}\bigtriangleup_{I;\kappa}^{\eta
}f\left(  x\right)  ,\ \ \ \ \ \text{for all }f\in L^{p},\label{bounded below}%
\\
\text{where }\bigtriangleup_{I;\kappa}^{\eta}f  &  \equiv\sum_{a\in\Gamma_{n}%
}\left\langle \left(  S_{\kappa,\eta}^{\mathcal{D}}\right)  ^{-1}%
f,h_{I;\kappa}^{a}\right\rangle \ h_{I;\kappa}^{a,\eta}\ ,\nonumber
\end{align}
and with convergence of the sum in both the $L^{p}$ norm and almost
everywhere, and%
\begin{equation}
\frac{1}{C_{p,n,\eta}}\left\Vert f\right\Vert _{L^{p}}\leq\left\Vert \left(
\sum_{I\in\mathcal{D}}\left\vert \bigtriangleup_{I;\kappa}^{\eta}f\right\vert
^{2}\right)  ^{\frac{1}{2}}\right\Vert _{L^{p}}\leq C_{p,n,\eta}\left\Vert
f\right\Vert _{L^{p}},\ \ \ \ \ \text{for all }f\in L^{p}.\label{square est}%
\end{equation}
Moreover, the smooth Alpert wavelets $\left\{  h_{I;\kappa}^{a,\eta}\right\}
_{I\in\mathcal{D},\ a\in\Gamma_{n}}$ are translation and dilation invariant in
the sense that $h_{I;\kappa}^{a,\eta}$ is a translate and dilate of the mother
Alpert wavelet $h_{I_{0};\kappa}^{a,\eta}$ where $I_{0}$ is the unit cube in
$\mathbb{R}^{n}$.
\end{theorem}

\begin{notation}
\label{Notation Alpert} We will often drop the index $a$ that parameterizes
the finite index set $\Gamma_{n}$ as it plays no essential role in most of
what follows, and it will be understood that when we write
\[
\bigtriangleup_{Q;\kappa}^{\eta}f=\left\langle \left(  S_{\kappa,\eta
}^{\mathcal{D}}\right)  ^{-1}f,h_{Q;\kappa}\right\rangle h_{Q;\kappa}^{\eta
}=\widehat{f}\left(  Q\right)  h_{Q;\kappa}^{\eta},
\]
we \emph{actually} mean the Alpert \emph{pseudoprojection},%
\[
\bigtriangleup_{Q;\kappa}^{\eta}f=\sum_{a\in\Gamma_{n}}\left\langle \left(
S_{\kappa,\eta}^{\mathcal{D}}\right)  ^{-1}f,h_{Q;\kappa}^{a}\right\rangle
h_{Q;\kappa}^{\eta,a}=\sum_{a\in\Gamma_{n}}\widehat{f_{a}}\left(  Q\right)
h_{Q;\kappa}^{a,\eta}\ ,
\]
where $\widehat{f_{a}}\left(  Q\right)  $ is a convenient abbreviation for the
inner product $\left\langle \left(  S_{\eta}^{\mathcal{D}}\right)
^{-1}f,h_{Q;\kappa}^{a}\right\rangle $ when $\kappa$ is understood. More
precisely, one can view $\widehat{f}\left(  Q\right)  =\left\{  \widehat
{f_{a}}\left(  Q\right)  \right\}  _{a\in\Gamma_{n}}$ and $h_{Q;\kappa}^{\eta
}=\left\{  h_{Q;\kappa}^{a,\eta}\right\}  _{a\in\Gamma_{n}}$ as finite
sequences of numbers and functions indexed by $\Gamma_{n}$, in which case
$\widehat{f}\left(  Q\right)  h_{Q;\kappa}^{\eta}$ is the dot product of these
two sequences. No confusion should arise between the Alpert coefficient
$\widehat{g}\left(  Q\right)  $, $Q\in\mathcal{G}\left[  U\right]  $ and the
Fourier transform $\widehat{g}\left(  \xi\right)  $, $\xi\in\mathbb{R}^{3}$,
as the argument in the first is a square in $\mathcal{G}\left[  U\right]  $,
while the argument in the second is a point in $\mathbb{R}^{3}$.
\end{notation}

\section{Equivalence of linear and trilinear single scale inequalities}

For $0<\delta<1$, let%
\begin{equation}
B_{\delta}\left(  0,2^{s}\right)  \equiv B\left(  0,2^{\frac{s}{1-\delta}%
}\right)  .\label{A delta}%
\end{equation}
Denote by $\mathcal{E}\left(  \otimes_{1}L^{q}\rightarrow L^{\frac{q}{3}%
}\right)  $ the Fourier extension inequality (\ref{FEC}).

\begin{definition}
Suppose $0<\delta<1$ and $\kappa\in\mathbb{N}$. Let $\varepsilon>0$,
$s\in\mathbb{N}$ and $1<q<\infty$. Denote by $\mathcal{A}^{s,\kappa,\delta
}\left(  \otimes_{1}L^{q}\rightarrow L^{q};\varepsilon\right)  $ the
\emph{local single scale} linear Fourier extension inequality%
\[
\left\Vert \mathcal{E}\mathsf{Q}_{s,U}^{\eta}f\right\Vert _{L^{q}\left(
B_{\delta}\left(  0,2^{s}\right)  \right)  }\leq C_{\delta,\varepsilon
,q}2^{\varepsilon s}\left\Vert f\right\Vert _{L^{q}\left(  U\right)  }\ ,
\]
taken over all single scale smooth Alpert pseudoprojections $\mathsf{Q}%
_{s,U}^{\eta}$ with $\kappa$ vanishing moments, and all $f\in L^{q}\left(
U\right)  $.
\end{definition}

Here is the corresponding \emph{local single scale} disjoint trilinear
inequality. Let $B_{\delta}\left(  0,2^{s}\right)  \equiv B\left(
0,2^{\frac{s}{1-\delta}}\right)  $ be the ball of radius $2^{\frac{s}%
{1-\delta}}$.

\begin{definition}
Suppose $0<\delta<1$ and $\kappa\in\mathbb{N}$. Let $\varepsilon>0$,
$s\in\mathbb{N}$ and $1<q<\infty$. Denote by $\mathcal{A}%
_{\operatorname*{disj}\nu}^{s,\kappa,\delta}\left(  \otimes_{3}L^{\infty
}\rightarrow L^{\frac{q}{3}};\varepsilon\right)  $ the local single scale
disjoint trilinear inequality%
\[
\left\Vert \mathcal{E}\mathsf{Q}_{s,U_{1}}^{\eta}f_{1}\ \mathcal{E}%
\mathsf{Q}_{s,U_{2}}^{\eta}f_{2}\ \mathcal{E}\mathsf{Q}_{s,U_{3}}^{\eta}%
f_{3}\right\Vert _{L^{\frac{q}{3}}\left(  B_{\delta}\left(  0,2^{s}\right)
\right)  }\leq C_{\delta,\kappa,\varepsilon,\nu}2^{\varepsilon s}\left\Vert
f_{1}\right\Vert _{L^{\infty}\left(  U_{1}\right)  }\left\Vert f_{2}%
\right\Vert _{L^{\infty}\left(  U_{2}\right)  }\left\Vert f_{3}\right\Vert
_{L^{\infty}\left(  U_{3}\right)  }\ ,
\]
taken over all $f_{k}\in L^{\infty}\left(  U_{k}\right)  $, all smooth Alpert
pseudoprojections $\mathsf{Q}_{s,U_{k}}^{\eta}$ with $\kappa$ vanishing
moments, and all $\nu$-disjoint triples $\left(  U_{1},U_{2},U_{3}\right)  \in
U^{3}$.
\end{definition}

Note that in this inequality we are testing the extension operator over smooth
Alpert pseudoprojections of bounded functions $f_{k}$, The following theorem
is the main technical result of this paper.

\begin{theorem}
[single scale Fourier and disjoint trilinear Fourier]\label{without}Suppose
$0<\delta<1$ and $\kappa\in\mathbb{N}$ with $\kappa>\frac{20}{\delta}$. Then
the local single scale Fourier extension inequality $\mathcal{A}%
_{\operatorname*{disj}\nu}^{s,\kappa,\delta}\left(  \otimes_{1}L^{\infty
}\rightarrow L^{q};\varepsilon\right)  $ for the paraboloid $\mathbb{P}^{2}$
in $\mathbb{R}^{3}$ holds \emph{if and only if} for every $q>3$ there is
$\nu>0$ depending on $q$, such that the local single scale $\nu$%
-disjoint\emph{\ }trilinear inequality $\mathcal{A}_{\operatorname*{disj}\nu
}^{s,\kappa,\delta}\left(  \otimes_{3}L^{\infty}\rightarrow L^{\frac{q}{3}%
};\varepsilon\right)  $ holds for all $s\in\mathbb{N}$ and $\varepsilon>0$.
\end{theorem}

The proof of Theorem \ref{without}\ is a modification of the corresponding
proof in \cite[proof of Theorem 3]{RiSa} (which was in turn based on the
simplest part of the clever argument in \cite[Section 2]{BoGu}) in which
bounded functions are now replaced with smooth Alpert pseudoprojections at a
single scale. More precisely, for $1<q\leq p<\infty$, $s\in\mathbb{N}$ and
$0<\delta<1$, we replace the quantity
\begin{equation}
Q_{R}^{\left(  q\right)  }\equiv\sup_{\left\Vert f\right\Vert _{L^{\infty
}\left(  U\right)  }\leq1}\left(  \int_{B\left(  0,R\right)  }\left\vert
\widehat{\Phi_{\ast}f}\left(  \xi\right)  \right\vert ^{q}d\xi\right)
^{\frac{1}{q}}\label{def Q}%
\end{equation}
defined in \cite[proof of Theorem 3]{RiSa}, and originally in \cite[Section
2]{BoGu}, with the Alpert projection quantity
\[
P_{\kappa,\delta}^{\left(  q\right)  }\left(  s\right)  \equiv\sup
_{\mathcal{G}\in\mathbb{G}}\sup_{\left\Vert f\right\Vert _{L^{\infty}\left(
\sigma\right)  }\leq1}\left(  \int_{B_{\delta}\left(  0,2^{s}\right)
}\left\vert \widehat{\Phi_{\ast}\mathsf{Q}_{\kappa,s,U}^{\eta,\mathcal{G}}%
f}\left(  \xi\right)  \right\vert ^{q}d\xi\right)  ^{\frac{1}{q}},
\]
where we have written $\mathsf{Q}_{\kappa,s,U}^{\eta,\mathcal{G}}$ in place of
$\mathsf{Q}_{s,U}^{\eta}$ when we wish to emphasize the grid $\mathcal{G}$ and
vanishing moment parameter $\kappa$ used in the construction of smooth Alpert
wavelets in\ Theorem \ref{reproducing}. We will show that for $\delta>0$ we
have $P_{\kappa,\delta}^{\left(  q\right)  }\left(  s\right)  \lesssim
C_{q,\delta}2^{C_{d}\delta s}$ provided $\kappa>\frac{20}{\delta}$, $q>3$, and
the $\nu$-disjoint\emph{\ }trilinear inequality $\mathcal{A}%
_{\operatorname*{disj}\nu}^{s,\kappa,\delta}\left(  \otimes_{3}L^{\infty
}\rightarrow L^{\frac{q}{3}};\varepsilon\right)  $ holds for all
$s\in\mathbb{N}$ and $\varepsilon>0$. Here the pseudoprojection $\mathsf{Q}%
_{\kappa,s,U}^{\eta,\mathcal{G}}$ is taken with respect to smooth Alpert
wavelets in the grid $\mathcal{G}$ having moment vanishing parameter $\kappa$.
We will often drop the superscript $\mathcal{G}$ when the actual grid is
unimportant, but it will play a key role in the parabolic dilation invariance
described below.

First note that the $L^{\infty}$ norm of the pseudoprojection $\mathsf{Q}%
_{\kappa,s,U}^{\eta,\mathcal{G}}f$ is controlled by the $L^{p}$ norm times a
harmless constant for $p$ large. More precisely, for any $1<r<\infty$ there is
a positive constant $C_{r}$ such that%
\begin{equation}
\left\Vert \mathsf{Q}_{\kappa,s,U}^{\eta,\mathcal{G}}f\right\Vert _{L^{\infty
}}\leq C_{r}2^{\frac{2}{r}s}\left\Vert f\right\Vert _{L^{p}}\leq C_{r}%
2^{\frac{2}{r}s}\left\Vert f\right\Vert _{L^{\infty}},\label{penalty}%
\end{equation}
since%
\begin{align*}
\left\vert \widehat{f}\left(  I\right)  \right\vert  &  =\left\vert
\left\langle \left(  S_{\kappa,\eta}^{\mathcal{G}}\right)  ^{-1}f,h_{I;\kappa
}\right\rangle \right\vert \leq\left\Vert \left(  S_{\kappa,\eta}%
^{\mathcal{G}}\right)  ^{-1}f\right\Vert _{L^{r}}\left\Vert h_{I;\kappa
}\right\Vert _{L^{r^{\prime}}}\\
&  \lesssim\left\Vert f\right\Vert _{L^{r}\left(  U\right)  }\left\Vert
h_{I;\kappa}\right\Vert _{L^{\infty}}\left\Vert \mathbf{1}_{I}\right\Vert
_{L^{r^{\prime}}}\lesssim\left\Vert f\right\Vert _{L^{r}\left(  U\right)
}\ell\left(  I\right)  ^{-1}\ell\left(  I\right)  ^{\frac{2}{r^{\prime}}}%
=\ell\left(  I\right)  ^{\frac{1}{r^{\prime}}-\frac{1}{r}}\left\Vert
f\right\Vert _{L^{r}}\ ,
\end{align*}
which gives%
\[
\left\vert \bigtriangleup_{I;\kappa}^{\eta}f\right\vert =\left\vert
\widehat{f}\left(  I\right)  h_{I;\kappa}^{\eta}\right\vert \lesssim
\ell\left(  I\right)  ^{-\left(  \frac{1}{r}-\frac{1}{r^{\prime}}\right)
}\left\Vert f\right\Vert _{L^{r}}\ell\left(  I\right)  ^{-1}\mathbf{1}%
_{\left(  1+\eta\ell\left(  I\right)  \right)  I}=\ell\left(  I\right)
^{-\frac{2}{r}}\left\Vert f\right\Vert _{L^{r}}\mathbf{1}_{\left(  1+\eta
\ell\left(  I\right)  \right)  I},
\]
for $1<r<\infty$. We will start the proof essentially as in \cite[proof of
Theorem 3]{RiSa} with a bounded function $f$, which followed \cite[Section
2]{BoGu} verbatim at this point, and then proceed with \textbf{Case 1} in
essentially the same way as well, but with an additional penalty factor
$2^{\frac{6}{qr}}$ for $r$ as large as we wish, due to the fact that we only
know that the map $\left(  S_{\kappa,\eta}^{\mathcal{G}}\right)  ^{-1}$ is
bounded on $L^{r}$ for $r<\infty$.

However, when using the quantity $P_{\kappa,\delta}^{\left(  q\right)
}\left(  s\right)  $ in place of $Q_{R}^{\left(  q\right)  }$,\ we must pay
extra attention to rescaling. The key to using the rescaling identity of Tao,
Vargas and Vega in \cite{TaVaVe} (that is used in the \textbf{Case 2} argument
below) is that the smooth Alpert pseudoprojection at a given scale, when
rescaled by a dyadic number, is again a smooth Alpert pseudoprojection in
another grid at scale corresponding to the rescaled ball $B_{\delta}\left(
0,2^{s}\right)  $. For this we will use $L^{p}$ normalized rescalings in order
to prove the case $p=\infty$, and then this is a consequence of the
translation and dilation invariance of the smooth Alpert wavelets
$h_{I;\kappa}^{\eta}$, and of a dilation invariance property of the maps
$S_{\kappa,\eta}^{\mathcal{G}}$.

The main difficulty here is the appearance of the operator $\left(
S_{\kappa,\eta}^{\mathcal{G}}\right)  ^{-1}$ in the inner product in the
projection $\bigtriangleup_{I;\kappa}^{\eta}f=\left\langle \left(
S_{\kappa,\eta}^{\mathcal{G}}\right)  ^{-1}f,h_{I;\kappa}\right\rangle
h_{I;\kappa}^{\eta}$. However, the following dilation invariance property of
$\left(  S_{\kappa,\eta}^{\mathcal{G}}\right)  ^{-1}$ will overcome this
difficulty. Let $\delta_{m,p}f\left(  x\right)  \equiv\frac{1}{2^{\frac{2}%
{p}m}}f\left(  \frac{x}{2^{m}}\right)  $ be the $L^{p}$ preserving dilation in
$\mathbb{R}^{2}$. Given a grid $\mathcal{G}$, let $2^{m}\mathcal{G}%
\equiv\left\{  2^{m}I:I\in\mathcal{G}\right\}  $ be the grid obtained by
dilating the squares in $\mathcal{G}$ by a factor $2^{m}$ about a fixed point
in $U$, that will vary from time to time, but is always understood from context.

\begin{lemma}
Let $1<q<p<\infty$. Then%
\begin{equation}
\delta_{m,p}\left(  S_{\kappa,\eta}^{\mathcal{G}}\right)  ^{-1}=\left(
S_{\kappa,\eta}^{2^{m}\mathcal{G}}\right)  ^{-1}\delta_{m,p}\ .\label{claim}%
\end{equation}

\end{lemma}

\begin{proof}
Since $f=\sum_{I\in\mathcal{G}\left[  U\right]  }\bigtriangleup_{I;\kappa
}^{\eta}f$ in $L^{p}$ by Theorem \ref{reproducing} when $p<\infty$, it
suffices to prove that%
\begin{equation}
\delta_{m,p}\left(  S_{\kappa,\eta}^{\mathcal{G}}\right)  ^{-1}\bigtriangleup
_{I;\kappa}^{\eta}f=\left(  S_{\kappa,\eta}^{2^{m}\mathcal{G}}\right)
^{-1}\delta_{m,p}\bigtriangleup_{I;\kappa}^{\eta}f\ ,\label{suff}%
\end{equation}
for all $I\in\mathcal{G}\left[  U\right]  $ and $f\in L^{p}\left(  U\right)
$. We begin by noting that the left and right hand sides of (\ref{suff})
satisfy
\begin{align*}
&  \delta_{m,p}\left(  S_{\kappa,\eta}^{\mathcal{G}}\right)  ^{-1}%
\bigtriangleup_{I;\kappa}^{\eta}f=\delta_{m,p}\left(  S_{\kappa,\eta
}^{\mathcal{G}}\right)  ^{-1}\left\langle \left(  S_{\kappa,\eta}%
^{\mathcal{G}}\right)  ^{-1}f,h_{I;\kappa}\right\rangle h_{I;\kappa}^{\eta}\\
&  =\left\langle \left(  S_{\kappa,\eta}^{\mathcal{G}}\right)  ^{-1}%
f,h_{I;\kappa}\right\rangle \delta_{m,p}\left(  S_{\kappa,\eta}^{\mathcal{G}%
}\right)  ^{-1}h_{I;\kappa}^{\eta}=\left\langle \left(  S_{\kappa,\eta
}^{\mathcal{G}}\right)  ^{-1}f,h_{I;\kappa}\right\rangle \delta_{m,p}%
h_{I;\kappa},
\end{align*}
and%
\begin{align*}
&  \left(  S_{\kappa,\eta}^{2^{m}\mathcal{G}}\right)  ^{-1}\delta
_{m,p}\bigtriangleup_{I;\kappa}^{\eta}f=\left(  S_{\kappa,\eta}^{2^{m}%
\mathcal{G}}\right)  ^{-1}\delta_{m,p}\left\langle \left(  S_{\kappa,\eta
}^{\mathcal{G}}\right)  ^{-1}f,h_{I;\kappa}\right\rangle h_{I;\kappa}^{\eta}\\
&  =\left\langle \left(  S_{\kappa,\eta}^{\mathcal{G}}\right)  ^{-1}%
f,h_{I;\kappa}\right\rangle \left(  S_{\kappa,\eta}^{2^{m}\mathcal{G}}\right)
^{-1}\delta_{m,p}h_{I;\kappa}^{\eta},
\end{align*}
respectively, and hence coincide once we have shown that%
\begin{equation}
\left(  S_{\kappa,\eta}^{2^{m}\mathcal{G}}\right)  ^{-1}\delta_{m,p}%
h_{I;\kappa}^{\eta}-\delta_{m,p}h_{I;\kappa}=0.\label{calc'}%
\end{equation}

To see (\ref{calc'}), we have%
\begin{align*}
\left(  S_{\kappa,\eta}^{2^{m}\mathcal{G}}\right)  ^{-1}\delta_{m,p}%
h_{I;\kappa}^{\eta}  &  =\left(  S_{\kappa,\eta}^{2^{m}\mathcal{G}}\right)
^{-1}\frac{1}{2^{\frac{2}{p}m}}h_{I;\kappa}^{\eta}\left(  \frac{x}{2^{m}%
}\right) \\
&  =\left(  S_{\kappa,\eta}^{2^{m}\mathcal{G}}\right)  ^{-1}\frac{1}%
{2^{\frac{2}{p}m}}2^{m}h_{2^{m}I;\kappa}^{\eta}\left(  x\right)  =2^{m\left(
1-\frac{2}{p}\right)  }h_{2^{m}I;\kappa}\left(  x\right)  ,
\end{align*}
since
\[
h_{I;\kappa}^{\eta}\left(  \frac{x}{2^{m}}\right)  =\frac{\sqrt{\left\vert
2^{m}I\right\vert }}{\sqrt{\left\vert I\right\vert }}h_{2^{m}I;\kappa}^{\eta
}\left(  x\right)  =2^{m}h_{2^{m}I;\kappa}^{\eta}\left(  x\right)  ,
\]
by the dilation and translation invariance of the smooth Alpert wavelets in
Theorem \ref{reproducing}, and where we have used that $2^{m}\mathcal{G}$ is
again a grid (not necessarily the same as $\mathcal{G}$), and that $2^{m}%
I\in2^{m}\mathcal{G}$. Hence the left hand side of (\ref{calc'}) equals
\begin{align*}
&  2^{m\left(  1-\frac{2}{p}\right)  }h_{2^{m}I;\kappa}\left(  x\right)
-\delta_{m,p}h_{I;\kappa}=2^{m\left(  1-\frac{2}{p}\right)  }h_{2^{m}I;\kappa
}\left(  x\right)  -2^{-\frac{2}{p}m}h_{I;\kappa}\left(  \frac{x}{2^{m}%
}\right) \\
&  =2^{m\left(  1-\frac{2}{p}\right)  }h_{2^{m}I;\kappa}\left(  x\right)
-2^{-\frac{2}{p}m}\frac{\sqrt{\left\vert 2^{m}I\right\vert }}{\sqrt{\left\vert
I\right\vert }}h_{2^{m}I;\kappa}\left(  x\right)  =0,
\end{align*}
which proves (\ref{calc'}).
\end{proof}

The proof of Theorem \ref{without} will use the dilation invariance property
(\ref{claim}), and is given in detail in the next section below. But first we
end the present section by proving the main Theorem \ref{single main Alpert}
assuming the key Theorem \ref{without}.

\subsection{Reduction to the characterization in \cite{RiSa}}

Here we will obtain the trilinear characterization of the Fourier extension
conjecture, as given in Theorem \ref{single main Alpert} of the introduction,
by using Theorem \ref{without} to reduce it to the multiscale disjoint
trilinear characterization in \cite{RiSa}, which we now recall.

\begin{definition}
Let $\varepsilon,\nu>0$ and $0<\delta<1$, $\kappa\in\mathbb{N}$ and
$1<q<\infty$. Denote by $\mathcal{A}_{\operatorname*{disj}\nu}%
^{\operatorname*{multi},\kappa,\delta}\left(  \otimes_{3}L^{\infty}\rightarrow
L^{\frac{q}{3}};\varepsilon\right)  $ the \emph{local multiscale} $\nu
$-disjoint trilinear smooth Alpert inequality,%
\begin{equation}
\left(  \int_{A\left(  0,2^{r}\right)  }\left(  \left\vert \mathcal{E}%
\mathsf{Q}_{\kappa,s_{1},U_{1}}^{\eta}f_{1}\left(  \xi\right)  \right\vert
\ \left\vert \mathcal{E}\mathsf{Q}_{\kappa,s_{2},U_{2}}^{\eta}f_{2}\left(
\xi\right)  \right\vert \ \left\vert \mathcal{E}\mathsf{Q}_{\kappa,s_{3}%
,U_{3}}^{\eta}f_{3}\left(  \xi\right)  \right\vert \right)  ^{\frac{q}{3}%
}\ d\xi\right)  ^{\frac{3}{q}}\leq\left(  S_{\varepsilon}^{\left(  q\right)
}2^{\varepsilon r}\right)  ^{3}\left\Vert f_{1}\right\Vert _{L^{\infty}%
}\left\Vert f_{2}\right\Vert _{L^{\infty}}\left\Vert f_{3}\right\Vert
_{L^{\infty}}\ ,\label{annular}%
\end{equation}
for all $r\in\mathbb{N}$, all $\nu$\emph{-disjoint} triples $\left(
U_{1},U_{2},U_{3}\right)  $, all smooth Alpert pseudoprojections
$\mathsf{Q}_{\kappa,s_{k},U_{k}}^{\eta}$ with moment parameter $\kappa$ and
with $s_{1}\leq s_{2}$ and $\frac{r}{1+\delta}<s_{2}\leq s_{3}<\frac
{r}{1-\delta} $, and all $f_{1},f_{2},f_{3}\in L^{\infty}$. Here
$S_{\varepsilon}^{\left(  q\right)  }$ denotes the best constant in the
inquality (\ref{annular}).
\end{definition}

\begin{theorem}
[\cite{RiSa}]\label{main Alpert}Let $0<\delta<1$ and $\kappa>\frac{20}{\delta
}$. The Fourier extension conjecture (\ref{FEC}) for the paraboloid in
$\mathbb{R}^{3}$ holds \emph{if and only if} for every $q>3$ there is $\nu>0$
depending only on $q$, such that $\mathcal{A}_{\operatorname*{disj}\nu
}^{\operatorname*{multi},\kappa,\delta}\left(  \otimes_{3}L^{\infty
}\rightarrow L^{\frac{q}{3}};\varepsilon\right)  $ holds for all
$\varepsilon>0$.
\end{theorem}

However, we can extract more\ from the proof given in \cite{RiSa}, namely a
quantitative inequality bounding the local norm $Q_{R}^{\left(  q\right)  }$
defined in (\ref{def Q}), by the local norm $S_{\varepsilon}^{\left(
q\right)  }2^{\varepsilon r}$ in (\ref{annular}), namely that there is a
positive constant $C$ such that%
\begin{equation}
Q_{2^{r}}^{\left(  q\right)  }\lesssim S_{\varepsilon}^{\left(  q\right)
}2^{C\varepsilon r},\ \ \ \ \ \text{for all }r\geq0\text{ and }%
q>3.\label{Q=<S}%
\end{equation}
We leave this routine and standard exercise for the reader. With Theorems
\ref{without} and \ref{main Alpert} in hand, we are now equipped to prove
Theorem \ref{single main Alpert}.

\subsubsection{Proof of Theorem \ref{single main Alpert}}

The Fourier extension conjecture (\ref{FEC}) trivially implies that the\ local
single scale Fourier extension inequality $\mathcal{A}_{\operatorname*{disj}%
\nu}^{s,\kappa,\delta}\left(  \otimes_{1}L^{\infty}\rightarrow L^{q}%
;\varepsilon\right)  $ holds for all $s\in\mathbb{N}$ and $\varepsilon>0$,
which by the easy direction of Theorem \ref{without}, implies that the local
single scale $\nu$-disjoint\emph{\ }trilinear inequality $\mathcal{A}%
_{\operatorname*{disj}\nu}^{s,\kappa,\delta}\left(  \otimes_{3}L^{\infty
}\rightarrow L^{\frac{q}{3}};\varepsilon\right)  $ holds for all
$s\in\mathbb{N}$ and $\varepsilon>0$. This proves the `only if' assertion of
Theorem \ref{single main Alpert}.

Conversely, assume the local single scale $\nu$-disjoint\emph{\ }trilinear
inequality $\mathcal{A}_{\operatorname*{disj}\nu}^{s,\kappa,\delta}\left(
\otimes_{3}L^{\infty}\rightarrow L^{\frac{q}{3}};\varepsilon\right)  $ holds
for all $s\in\mathbb{N}$ and $\varepsilon>0$. Then the hard direction of
Theorem \ref{without} shows that the local single scale linear Fourier
extension inequality $\mathcal{A}_{\operatorname*{disj}\nu}^{s,\kappa,\delta
}\left(  \otimes_{1}L^{\infty}\rightarrow L^{q};\varepsilon\right)  $ holds
for all $s\in\mathbb{N}$ and $\varepsilon>0$.

From $\mathcal{A}_{\operatorname*{disj}\nu}^{s,\kappa,\delta}\left(
\otimes_{1}L^{\infty}\rightarrow L^{q};\varepsilon\right)  $ and H\"{o}lder's
inequality, we now obtain
\begin{align*}
&  \int_{A\left(  0,2^{r}\right)  }\left(  \left\vert \mathcal{E}%
\mathsf{Q}_{\kappa,s_{1},U_{1}}^{\eta}f_{1}\left(  \xi\right)  \right\vert
\ \left\vert \mathcal{E}\mathsf{Q}_{\kappa,s_{2},U_{2}}^{\eta}f_{2}\left(
\xi\right)  \right\vert \ \left\vert \mathcal{E}\mathsf{Q}_{\kappa,s_{3}%
,U_{3}}^{\eta}f_{3}\left(  \xi\right)  \right\vert \right)  ^{\frac{q}{3}%
}\ d\xi\\
&  \leq\left(  \int_{A\left(  0,2^{r}\right)  }\left\vert \mathcal{E}%
\mathsf{Q}_{\kappa,s_{1},U_{1}}^{\eta}f_{k}\left(  \xi\right)  \right\vert
^{q}d\xi\right)  ^{\frac{1}{3}}\prod_{k=2}^{3}\left(  \int_{A\left(
0,2^{r}\right)  }\left\vert \mathcal{E}\mathsf{Q}_{\kappa,s_{k},U_{k}}^{\eta
}f_{k}\left(  \xi\right)  \right\vert ^{q}d\xi\right)  ^{\frac{1}{3}}\\
&  \lesssim\left(  Q_{2^{r}}^{\left(  q\right)  }\left\Vert f_{1}\right\Vert
_{L^{\infty}}\right)  ^{\frac{q}{3}}\left(  2^{\varepsilon r}\left\Vert
f_{2}\right\Vert _{L^{\infty}}2^{\varepsilon r}\left\Vert f_{3}\right\Vert
_{L^{\infty}\left(  U_{3}\right)  }\right)  ^{\frac{q}{3}}\\
&  \lesssim\left(  S_{\varepsilon}^{\left(  q\right)  }2^{C\varepsilon
r}\right)  ^{\frac{q}{3}}\left\Vert f_{1}\right\Vert _{L^{\infty}}^{\frac
{q}{3}}\left(  2^{\varepsilon r}\left\Vert f_{2}\right\Vert _{L^{\infty}%
}2^{\varepsilon r}\left\Vert f_{3}\right\Vert _{L^{\infty}\left(
U_{3}\right)  }\right)  ^{\frac{q}{3}},
\end{align*}
where we have used the definition of $Q_{2^{r}}^{\left(  q\right)  }$ in
(\ref{def Q}) in the third line, and then (\ref{Q=<S}) in the final line. Now
take the $q^{th}$ root of the above inequality, and then take the supremum
over all $s_{1}\leq s_{2}$ and $\frac{r}{1+\delta}<s_{2}\leq s_{3}<\frac
{r}{1-\delta}$, and all $f_{k}\in L^{\infty}\left(  U_{k}\right)  $ with norm
$1$, to conclude that%
\[
S_{\varepsilon}^{\left(  q\right)  }=\sup_{s_{k,}f_{k}}\left(  \int_{A\left(
0,2^{r}\right)  }\left(  \left\vert \mathcal{E}\mathsf{Q}_{\kappa,s_{1},U_{1}%
}^{\eta}f_{1}\left(  \xi\right)  \right\vert \ \left\vert \mathcal{E}%
\mathsf{Q}_{\kappa,s_{2},U_{2}}^{\eta}f_{2}\left(  \xi\right)  \right\vert
\ \left\vert \mathcal{E}\mathsf{Q}_{\kappa,s_{3},U_{3}}^{\eta}f_{3}\left(
\xi\right)  \right\vert \right)  ^{\frac{q}{3}}\ d\xi\right)  ^{\frac{1}{q}%
}\lesssim\left(  S_{\varepsilon}^{\left(  q\right)  }2^{C\varepsilon
r}\right)  ^{\frac{1}{3}}\left(  2^{\varepsilon r}\right)  ^{\frac{2}{3}}.
\]
Since $S_{\varepsilon}^{\left(  q\right)  }<\infty$, we can divide through by
$\left(  S_{\varepsilon}^{\left(  q\right)  }\right)  ^{\frac{1}{3}}$ to
obtain $\left(  S_{\varepsilon}^{\left(  q\right)  }\right)  ^{\frac{2}{3}%
}\lesssim2^{\frac{C+2}{3}\varepsilon r}$, and then another application of
(\ref{Q=<S}) yields%
\[
Q_{2^{r}}^{\left(  q\right)  }\lesssim S_{\varepsilon}^{\left(  q\right)
}2^{C\varepsilon r}\lesssim2^{q\frac{C+2}{2}\varepsilon r}2^{C\varepsilon r},
\]
which is well known to imply the Fourier extension conjecture by factorization
and Tao's $\varepsilon$-removal, see e.g. \cite{BoGu}.

\section{Proof of Theorem \ref{without} on local single scale inequalities}

The proof of Theorem \ref{without} is more than just a simple reprise of the
argument in \cite{RiSa} adapted from \cite{BoGu}, due to the details
surrounding parabolic rescalings of smooth Alpert wavelets. Suppose $S$ is a
compact smooth hypersurface contained in $\mathbb{R}^{3}$ that is contained in
the paraboloid $\mathbb{P}^{2}$, and denote surface measure on $S$ by $\sigma
$. The next definition is specialized from \cite{BoGu}, but with their local
quantity $Q_{R}^{\left(  q\right)  }\equiv\sup_{\left\Vert f\right\Vert
_{L^{\infty}\left(  \sigma\right)  }\leq1}\left(  \int_{B\left(  0,R\right)
}\left\vert \widehat{f\sigma}\left(  \xi\right)  \right\vert ^{q}d\xi\right)
^{\frac{1}{q}}$ replaced by a local single scale Alpert quantity
$P_{\kappa,\delta}^{\left(  q\right)  }\left(  s\right)  $, in which the scale
$s$ of the pseudoprojection $\mathsf{Q}_{\kappa,s,U}^{\eta,\mathcal{G}}f$ is
tied to the size of the local ball%
\[
B_{\delta}\left(  0,2^{s}\right)  \equiv B\left(  0,2^{\frac{s}{1-\delta}%
}\right)  .
\]
Thus $2^{s}$ plays the role of $R$ here. Recall also that we write
\[
\mathsf{Q}_{s,U}^{\eta}f=\mathsf{Q}_{\kappa,s,U}^{\eta,\mathcal{G}}f\equiv
\sum_{I\in\mathcal{G}_{s}\left[  U\right]  }\bigtriangleup_{I;\kappa}^{\eta
}f=\sum_{I\in\mathcal{G}_{s}\left[  U\right]  }\left\langle \left(
S_{\kappa,\eta}\right)  ^{-1}f,h_{I;\kappa}\right\rangle h_{I;\kappa}^{\eta
}\ ,
\]
when we wish to emphasize the grid $\mathcal{G}$ and vanishing moment
parameter $\kappa$ used in the construction of smooth Alpert wavelets
in\ Theorem \ref{reproducing}. Here $\left\{  h_{I;\kappa}\right\}
_{I\in\mathcal{G}}$ is the orthonormal basis of Alpert wavelets $h_{I;\kappa}$
in $L^{2}\left(  U\right)  $, and $\left\{  h_{I;\kappa}^{\eta}\right\}
_{I\in\mathcal{G}}$ is the corresponding frame of smooth Alpert wavelets
$h_{I;\kappa}^{\eta}$ with $\kappa$ vanishing moments in $L^{p}\left(
U\right)  $, $1<p<\infty$.

\begin{definition}
For $0<\delta<1$, $1<q<\infty$ and $\kappa,s\in\mathbb{N}$ define
$P_{\kappa,\delta}^{\left(  q\right)  }\left(  s\right)  $ to be the best
constant in the local linear Fourier extension inequality,%
\[
\left(  \int_{B_{\delta}\left(  0,2^{s}\right)  }\left\vert \mathcal{E}\left(
\Phi_{\ast}\mathsf{Q}_{\kappa,s,U}^{\eta,\mathcal{G}}f\right)  \left(
\xi\right)  \right\vert ^{q}d\xi\right)  ^{\frac{1}{q}}\leq P_{\kappa,\delta
}^{\left(  q\right)  }\left(  s\right)  \left\Vert f\right\Vert _{L^{\infty
}\left(  U\right)  }\ ,
\]
taken over all grids $\mathcal{G}\in\mathbb{G}$ and$\ $all $f\in L^{\infty
}\left(  U\right)  $, i.e.
\begin{equation}
P_{\kappa,\delta}^{\left(  q\right)  }\left(  s\right)  \equiv\sup
_{\mathcal{G}\in\mathbb{G}}\sup_{\left\Vert f\right\Vert _{L^{\infty}\left(
U\right)  }\leq1}\left(  \int_{B_{\delta}\left(  0,2^{s}\right)  }\left\vert
\mathcal{E}\left(  \Phi_{\ast}\mathsf{Q}_{\kappa,s,U}^{\eta,\mathcal{G}%
}f\right)  \left(  \xi\right)  \right\vert ^{q}d\xi\right)  ^{\frac{1}{q}%
}.\label{Q_R}%
\end{equation}

\end{definition}

Here is the main step in the proof of Theorem \ref{without}.

\begin{theorem}
\label{Loc lin}Let $\varepsilon,\nu>0$, $0<\delta<1$, $\kappa\in\mathbb{N}$,
$1<q<\infty$, and let $S$ and $U$ be as above. Then for every $\varepsilon>0$
and $q>3$, there is $\nu=\nu\left(  q\right)  >0$ depending only on $q>3$ and
a universal constant $C_{0}$ such that if $\mathcal{A}_{\operatorname*{disj}%
\nu}^{s,\kappa,\delta}\left(  \otimes_{3}L^{\infty}\rightarrow L^{\frac{q}{3}%
};\varepsilon\right)  $ holds for all $s\in\mathbb{N}$ and $\varepsilon>0$,
then%
\[
P_{\kappa,\delta}^{\left(  q\right)  }\left(  s\right)  \leq C_{\kappa
,\delta,q}+C_{\kappa,q,\delta,\varepsilon}2^{\varepsilon s}%
,\ \ \ \ \ \text{for all }s>C_{0}\frac{q}{q-3}\text{ and }\varepsilon>0.
\]

\end{theorem}

\subsection{Dilation invariance}

To prove Theorem \ref{Loc lin}, we will need to establish a dilation
invariance property of smooth Alpert pseudoprojections of $L^{\infty}$
functions. Since we only know that the operators $\left(  S_{\kappa,\eta
}^{\mathcal{G}}\right)  ^{-1}$ are bounded on $L^{p}$ with $p<\infty$, it will
be convenient to include such $p$ in our analysis.

Define%
\[
\mathcal{F}_{\kappa,p}\left(  s\right)  \equiv\left\{  \mathsf{Q}_{\kappa
,s,U}^{\eta,\mathcal{G}}f:\mathcal{G}\in\mathbb{G}\text{ and }\left\Vert
f\right\Vert _{L^{p}\left(  U\right)  }\leq1\right\}  ,
\]
so that%
\[
P_{\kappa,\delta}^{\left(  q\right)  }\left(  s\right)  =\sup_{F\in
\mathcal{F}_{\kappa,\infty}\left(  s\right)  }\left(  \int_{B_{\delta}\left(
0,2^{s}\right)  }\left\vert \mathcal{E}\left(  \Phi_{\ast}F\right)  \left(
\xi\right)  \right\vert ^{q}d\xi\right)  ^{\frac{1}{q}}.
\]
We now prove that the families $\mathcal{F}_{\kappa,p}\left(  s\right)  $,
$s\in\mathbb{N}$, satisfy a crucial invariance property under the dyadic
$L^{p}$ preserving dilations $\delta_{m,p}f\left(  x\right)  \equiv\frac
{1}{2^{\frac{2}{p}m}}f\left(  \frac{x}{2^{k}}\right)  $ in $\mathbb{R}^{2}$,
i.e. $\delta_{m,p}\mathcal{F}_{\kappa,p}\left(  s\right)  =\mathcal{F}%
_{\kappa,p}\left(  s-m\right)  $.

But first we prove the dilation invariance of smooth Alpert pseudoprojections.
Recall the dilated grid%
\[
2^{m}\mathcal{G}\equiv\left\{  2^{m}I:I\in\mathcal{G}\right\}  ,
\]
and write $\bigtriangleup_{I;\kappa}^{\eta,\mathcal{G}}f=\left\langle \left(
S_{\kappa,\eta}^{\mathcal{G}}\right)  ^{-1}f,h_{I;\kappa}\right\rangle
h_{I;\kappa}^{\eta}$ when we wish to emphasize the underlying grid
$\mathcal{G}$ in the pseudoprojection $\bigtriangleup_{I;\kappa}^{\eta}f$.

\begin{lemma}
\label{proj inv}Let $1<q<p\leq\infty$. Then%
\begin{equation}
\delta_{m,p}\left(  \bigtriangleup_{I;\kappa}^{\eta,\mathcal{G}}f\right)
=\bigtriangleup_{2^{m}I;\kappa}^{\eta,2^{m}\mathcal{G}}\delta_{m,p}%
f,\ \ \ \ \ \text{if both }I,2^{m}I\in\mathcal{G}\left[  U\right]  \text{ and
}f\in L^{p}\left(  U\right)  ,\label{inv}%
\end{equation}

\end{lemma}

\begin{proof}
It suffices to prove the case $p<\infty$ since $L^{\infty}\left(  U\right)
\subset L^{p}\left(  U\right)  $. We have,%
\begin{equation}
\delta_{m,p}\left(  \bigtriangleup_{I;\kappa}^{\eta,\mathcal{G}}f\right)
\left(  x\right)  =\delta_{m,p}\left(  \left\langle \left(  S_{\kappa,\eta
}^{\mathcal{G}}\right)  ^{-1}f,h_{I;\kappa}\right\rangle h_{I;\kappa}^{\eta
}\right)  \left(  x\right)  =\left\langle \left(  S_{\kappa,\eta}%
^{\mathcal{G}}\right)  ^{-1}f,h_{I;\kappa}\right\rangle \delta_{m,p}%
h_{I;\kappa}^{\eta}\left(  x\right)  ,\label{one}%
\end{equation}
and%
\begin{equation}
\bigtriangleup_{2^{m}I;\kappa}^{\eta,2^{m}\mathcal{G}}\delta_{m,p}f\left(
x\right)  =\left\langle \left(  S_{\kappa,\eta}^{2^{m}\mathcal{G}}\right)
^{-1}\delta_{m,p}f,h_{2^{m}I;\kappa}\right\rangle h_{2^{m}I;\kappa}^{\eta
}\left(  x\right)  =\left\langle \delta_{m,p}\left(  S_{\kappa,\eta
}^{\mathcal{G}}\right)  ^{-1}f,h_{2^{m}I;\kappa}\right\rangle h_{2^{m}%
I;\kappa}^{\eta}\left(  x\right)  ,\label{two}%
\end{equation}
since $2^{m}I\in2^{m}\mathcal{G}$, and since $\left(  S_{\kappa,\eta}%
^{2^{m}\mathcal{G}}\right)  ^{-1}\delta_{m,p}=\delta_{m,p}\left(
S_{\kappa,\eta}^{\mathcal{G}}\right)  ^{-1}$ by (\ref{claim}).

We first compute the factor $\left\langle \delta_{m,p}\left(  S_{\kappa,\eta
}^{\mathcal{G}}\right)  ^{-1}f,h_{2^{m}I;\kappa}\right\rangle $ in (\ref{two})
in terms of $I$. We have%
\begin{align}
&  \ \ \ \ \ \ \ \ \ \ \ \ \ \ \ \ \ \ \ \ \ \ \ \left\langle \delta
_{m,p}\left(  S_{\kappa,\eta}^{\mathcal{G}}\right)  ^{-1}f,h_{2^{m}I;\kappa
}\right\rangle =\int\frac{1}{2^{\frac{2}{p}m}}\left(  S_{\kappa,\eta
}^{\mathcal{G}}\right)  ^{-1}f\left(  \frac{y}{2^{m}}\right)  h_{2^{m}%
I;\kappa}\left(  y\right)  dy\label{three}\\
&  =\int\frac{1}{2^{\frac{2}{p}m}}\left(  S_{\kappa,\eta}^{\mathcal{G}%
}\right)  ^{-1}f\left(  z\right)  h_{2^{m}I;\kappa}\left(  2^{m}z\right)
2^{2m}dz=2^{2m}\int\left(  S_{\kappa,\eta}^{\mathcal{G}}\right)  ^{-1}f\left(
z\right)  \frac{1}{2^{\frac{2}{p}m}}\frac{h_{2^{m}I;\kappa}\left(
2^{m}z\right)  }{h_{I;\kappa}\left(  z\right)  }h_{I;\kappa}\left(  z\right)
dz\nonumber\\
&  =2^{2m}\int\left(  S_{\kappa,\eta}^{\mathcal{G}}\right)  ^{-1}f\left(
z\right)  \frac{1}{2^{\frac{2}{p}m}}\sqrt{\frac{\left\vert I\right\vert
}{\left\vert 2^{m}I\right\vert }}h_{I;\kappa}\left(  z\right)  dz=\frac{2^{m}%
}{2^{\frac{2}{p}m}}\int\left(  S_{\kappa,\eta}^{\mathcal{G}}\right)
^{-1}f\left(  z\right)  h_{I;\kappa}\left(  z\right)  dz=2^{m\left(
1-\frac{2}{p}\right)  }\left\langle \left(  S_{\kappa,\eta}^{\mathcal{G}%
}\right)  ^{-1}f,h_{I;\kappa}\right\rangle .\nonumber
\end{align}
Next we compute the factor $\delta_{m,p}h_{I;\kappa}^{\eta}\left(  x\right)  $
in (\ref{one}) in terms of $2^{m}I$. We have%
\begin{equation}
\delta_{m,p}h_{I;\kappa}^{\eta}\left(  x\right)  =\frac{1}{2^{\frac{2}{p}m}%
}h_{I;\kappa}^{\eta}\left(  \frac{x}{2^{m}}\right)  =\frac{1}{2^{\frac{2}{p}%
m}}\frac{h_{I;\kappa}\left(  \frac{x}{2^{m}}\right)  }{h_{2^{m}I;\kappa
}\left(  x\right)  }h_{2^{m}I;\kappa}^{\eta}\left(  x\right)  =\frac
{1}{2^{\frac{2}{p}m}}\sqrt{\frac{\left\vert 2^{m}I\right\vert }{\left\vert
I\right\vert }}h_{2^{m}I;\kappa}^{\eta}\left(  x\right)  =2^{m\left(
1-\frac{2}{p}\right)  }h_{2^{m}I;\kappa}^{\eta}\left(  x\right)  .\label{four}%
\end{equation}
Using first (\ref{one}), then (\ref{four}), then (\ref{three}) and finally
(\ref{two}) gives (\ref{inv}),%
\begin{align*}
\delta_{m,p}\left(  \bigtriangleup_{I;\kappa}^{\eta,\mathcal{G}}f\right)
\left(  x\right)   &  =\left\langle \left(  S_{\kappa,\eta}^{\mathcal{G}%
}\right)  ^{-1}f,h_{I;\kappa}\right\rangle \delta_{m,p}h_{I;\kappa}^{\eta
}\left(  x\right) \\
&  =\left\langle \left(  S_{\kappa,\eta}^{\mathcal{G}}\right)  ^{-1}%
f,h_{I;\kappa}\right\rangle 2^{m\left(  1-\frac{2}{p}\right)  }h_{2^{m}%
I;\kappa}^{\eta}\left(  x\right) \\
&  =\left\langle \delta_{m,p}\left(  S_{\kappa,\eta}^{\mathcal{G}}\right)
^{-1}f,h_{2^{m}I;\kappa}\right\rangle h_{2^{m}I;\kappa}^{\eta}\left(  x\right)
\\
&  =\bigtriangleup_{2^{m}I;\kappa}^{\eta,2^{m}\mathcal{G}}\delta_{m,p}f\left(
x\right)  .
\end{align*}

\end{proof}

We immediately obtain the desired dilation invariance of\ the families
$\mathcal{F}_{\kappa,p}\left(  s\right)  $, $s\in\mathbb{N}$.

\begin{corollary}
\label{F inv}Let $1<q\leq p\leq\infty$. Then the sequence of families of
functions $\left\{  \mathcal{F}_{\kappa,p}\left(  s\right)  \right\}
_{s=1}^{\infty}$ is invariant under the dyadic $L^{p}$ preserving dilations
$\delta_{m,p}$ in the following sense,%
\[
\delta_{m,p}\mathcal{F}_{\kappa,p}\left(  s\right)  =\mathcal{F}_{\kappa
,p}\left(  s-m\right)  .
\]

\end{corollary}

\subsection{The pigeonholing argument of Bourgain and Guth\label{Sub pigeon}}

\begin{proof}
[Proof of Theorem \ref{Loc lin}]We begin the argument as in \cite{BoGu}, with
notation as in \cite{RiSa}. We let the surface $S$ be a compact smooth piece
of the paraboloid $\mathbb{P}^{2}$ given by $z_{3}=\left\vert z^{\prime
}\right\vert ^{2}=z_{1}^{2}+z_{2}^{2}$ in $\mathbb{R}^{3}$, and for $f\in
L^{\infty}\left(  S\right)  $ with $\left\Vert f\right\Vert _{L^{\infty
}\left(  S\right)  }=1$, we consider the oscillatory integral,%
\begin{align*}
Tf\left(  \xi\right)   &  =\int_{U}e^{i\phi\left(  \xi,y\right)  }f\left(
y\right)  dy=\int_{U}e^{i\left\{  \xi_{1}\cdot y_{1}+\xi_{2}\cdot y_{2}%
+\xi_{3}\left(  y_{1}^{2}+y_{2}^{2}\right)  \right\}  }f\left(  y\right)  dy\\
&  =\int_{U}e^{i\xi\cdot\left(  y,\left\vert y\right\vert ^{2}\right)
}f\left(  y\right)  dy=\widehat{f^{\Phi}}\left(  \xi\right)  =\widehat
{\Phi_{\ast}\left[  f\left(  y\right)  dy\right]  }\left(  \xi\right)
,\ \ \ \ \ \text{for }\xi\in\mathbb{R}^{3},
\end{align*}
where
\[
\phi\left(  \xi,y\right)  =\xi\cdot\Phi\left(  y\right)  \text{ and }%
\Phi\left(  y\right)  \equiv\left(  y_{1},y_{2},y_{1}^{2}+y_{2}^{2}\right)  .
\]
For $\lambda\geq1$, let $f=\sum_{I\in\mathcal{G}_{\lambda}\left[  U\right]
}\boldsymbol{1}_{I}f=\sum_{I\in\mathcal{G}_{\lambda}\left[  U\right]  }f_{I}$
and write%
\[
Tf\left(  \xi\right)  =\sum_{I\in\mathcal{G}_{\lambda}\left[  U\right]  }%
\int_{S}e^{i\xi\cdot\left(  y,\left\vert y\right\vert ^{2}\right)  }%
f_{I}\left(  y\right)  dy=\sum_{I\in\mathcal{G}_{\lambda}\left[  U\right]
}e^{i\phi\left(  \xi,c_{I}\right)  }\int e^{i\left\{  \phi\left(
\xi,y\right)  -\phi\left(  \xi,c_{I}\right)  \right\}  }f_{I}\left(  y\right)
dy=\sum_{I\in\mathcal{G}_{\lambda}\left[  S\right]  }e^{i\phi\left(  \xi
,c_{I}\right)  }T_{I}f\left(  \xi\right)  ,
\]
where
\begin{align}
T_{I}f\left(  \xi\right)   &  \equiv\int e^{i\left\{  \phi\left(
\xi,y\right)  -\phi\left(  \xi,c_{I}\right)  \right\}  }f_{I}\left(  y\right)
dy=e^{-i\left\{  \xi\cdot\Phi\left(  c_{I}\right)  \right\}  }\int
e^{i\xi\cdot\Phi\left(  y\right)  }f_{I}\left(  y\right)  dy\label{def T_I}\\
&  =e^{-i\left\{  \xi\cdot\Phi\left(  c_{I}\right)  \right\}  }\widehat
{f_{I}^{\Phi}}\left(  \xi\right)  =\widehat{\tau_{-\Phi\left(  c_{I}\right)
}f_{I}^{\Phi}}\left(  \xi\right)  ,\nonumber
\end{align}
where $\tau_{\Phi\left(  c_{I}\right)  }g\left(  z\right)  \equiv g\left(
z-\Phi\left(  c_{I}\right)  \right)  $ is translation of a function $g$ by the
unit vector $\Phi\left(  c_{I}\right)  $.

Note that%
\[
\left\vert \nabla_{\xi}\left\{  \xi\cdot\left(  y-c_{I},\left\vert
y\right\vert ^{2}-\left\vert c_{I}\right\vert ^{2}\right)  \right\}
\right\vert =\left\vert \left(  y-c_{I},\left\vert y\right\vert ^{2}%
-\left\vert c_{I}\right\vert ^{2}\right)  \right\vert \lesssim\frac
{1}{2^{\lambda}},\ \ \ \ \ \text{for }y\in I,
\]
implies%
\begin{align}
\nabla_{\xi}T_{I}f\left(  \xi\right)   &  =\nabla_{\xi}\int e^{i\xi
\cdot\left(  y-c_{I},\left\vert y-c_{I}\right\vert ^{2}\right)  }f_{I}\left(
y\right)  dy=\int\nabla_{\xi}e^{i\xi\cdot\left(  y-c_{I},\left\vert
y-c_{I}\right\vert ^{2}\right)  }f_{I}\left(  y\right)
dy\label{pre note that}\\
&  =\int ie^{i\xi\cdot\left(  y-c_{I},\left\vert y-c_{I}\right\vert
^{2}\right)  }\nabla_{\xi}\left\{  \xi\cdot\left(  y-c_{I},\left\vert
y\right\vert ^{2}-\left\vert c_{I}\right\vert ^{2}\right)  \right\}
f_{I}\left(  y\right)  dy,\nonumber
\end{align}
which implies%
\begin{equation}
\left\vert \nabla_{\xi}T_{I}f\left(  \xi\right)  \right\vert \leq
\int\left\vert \nabla_{\xi}\left\{  \xi\cdot\left(  y-c_{I},\left\vert
y\right\vert ^{2}-\left\vert c_{I}\right\vert ^{2}\right)  \right\}
\right\vert \left\vert f_{I}\left(  y\right)  \right\vert dy\lesssim\frac
{1}{2^{\lambda}}\left\Vert f_{I}\right\Vert _{L^{1}\left(  U\right)  }%
\lesssim\frac{1}{2^{3\lambda}}\left\Vert f_{I}\right\Vert _{L^{\infty}\left(
U\right)  },\label{note that}%
\end{equation}
since $\ell\left(  I\right)  =\frac{1}{2^{\lambda}}$. We will use the
estimates (\ref{pre note that}) and (\ref{note that}) in (\ref{exclaim'}) below.

Now let $\rho$ be a smooth rapidly decreasing bump function such that
$\widehat{\rho}\left(  \xi\right)  =1$ for $\left\vert \xi\right\vert \leq1$,
and set%
\[
\rho_{\lambda}\left(  z\right)  \equiv\frac{1}{2^{3\lambda}}\rho\left(
\frac{z}{2^{\lambda}}\right)  ,\ \ \ \ \ \widehat{\rho_{\lambda}}\left(
\xi\right)  =\widehat{\rho}\left(  2^{\lambda}\xi\right)  =1\text{ on
}B\left(  0,2^{-\lambda}\right)  \text{ and }\rho_{\lambda}\left(  z\right)
\approx\frac{1}{2^{3\lambda}}\text{ on }B\left(  0,2^{\lambda}\right)  .
\]
Then from (\ref{def T_I}) we obtain%
\[
T_{I}f\left(  \xi\right)  =T_{I}f\ast\rho_{\lambda}\left(  \xi\right)
\ ,\ \ \ \ \ \text{for }I\in\mathcal{G}_{\lambda}\left[  S\right]  \text{ and
}\xi\in\mathbb{R}^{3},
\]
since $\tau_{-\Phi\left(  c_{I}\right)  }f_{I}^{\Phi}\subset B\left(
0,2^{-\lambda}\right)  $ and $\widehat{T_{I}f}\left(  z\right)  =\tau
_{-\Phi\left(  c_{I}\right)  }f_{I}^{\Phi}\left(  z\right)  $ imply%
\[
\widehat{T_{I}f\ast\rho_{\lambda}}\left(  z\right)  =\widehat{T_{I}f}\left(
z\right)  \widehat{\rho_{\lambda}}\left(  z\right)  =\tau_{-\Phi\left(
c_{I}\right)  }f_{I}^{\Phi}\left(  z\right)  \widehat{\rho_{\lambda}}\left(
z\right)  =\tau_{-\Phi\left(  c_{I}\right)  }f_{I}^{\Phi}\left(  z\right)
=\widehat{T_{I}f}\left(  z\right)  .
\]

Now fix a point $a\in\Gamma_{\lambda}\left(  R\right)  $, where%
\[
\Gamma_{\lambda}\left(  R\right)  \equiv2^{\lambda}\mathbb{Z}^{3}\cap
B_{R},\ \ \ \ \ \text{and }B_{R}\equiv B\left(  0,R\right)  ,
\]
and restrict $\xi$ to the ball $B\left(  a,2^{\lambda}\right)  $. We will
write $B_{R}$ as apposed to $B\left(  0,R\right)  $ above in order to
emphasize the different roles played by $\lambda$ and $R$, namely
$R\nearrow\infty$ while $\lambda$ remains a fixed sufficiently large integer
to be chosen.

Then for $\xi\in B\left(  a,2^{\lambda}\right)  $ and $I\in\mathcal{G}%
_{\lambda}\left[  S\right]  $ we have
\begin{align}
&  \left\vert T_{I}f\left(  \xi\right)  \right\vert =\left\vert T_{I}f\ast
\rho_{\lambda}\left(  \xi\right)  \right\vert =\left\vert \int_{\mathbb{R}%
^{3}}T_{I}f\left(  z\right)  \rho_{\lambda}\left(  \xi-z\right)  dz\right\vert
\label{exclaim}\\
&  \leq\int_{\mathbb{R}^{3}}\left\vert T_{I}f\left(  z\right)  \right\vert
\left\vert \rho_{\lambda}\left(  \xi-z\right)  \right\vert dz\leq
\int_{\mathbb{R}^{3}}\left\vert T_{I}f\left(  z\right)  \right\vert
\sup_{\omega\in B\left(  a,2^{\lambda}\right)  }\left\vert \rho_{\lambda
}\left(  z-\omega\right)  \right\vert dz=\int_{\mathbb{R}^{3}}\left\vert
T_{I}f\left(  z\right)  \right\vert \zeta_{\lambda}\left(  z-a\right)
dz,\nonumber
\end{align}
where $\zeta_{\lambda}\left(  w\right)  \equiv\sup_{\omega-a\in B\left(
0,2^{\lambda}\right)  }\left\vert \rho_{\lambda}\left(  \omega\right)
\right\vert $, since $\rho$ can be chosen radial and,
\begin{align*}
\sup_{\omega\in B\left(  a,2^{\lambda}\right)  }\left\vert \rho_{\lambda
}\left(  \omega-z\right)  \right\vert  &  =\frac{1}{2^{3\lambda}}\sup
_{\omega\in B\left(  a,2^{\lambda}\right)  }\left\vert \rho\left(
\frac{\left(  \omega-a\right)  -\left(  z-a\right)  }{2^{\lambda}}\right)
\right\vert =\frac{1}{2^{3\lambda}}\sup_{\gamma\in B\left(  0,1\right)
}\left\vert \rho\left(  \frac{z-a}{2^{\lambda}}-\gamma\right)  \right\vert
=\zeta_{\lambda}\left(  z-a\right)  ,\\
\text{where }\zeta\left(  w\right)   &  \equiv\sup_{\left\vert w-w^{\prime
}\right\vert \leq1}\left\vert \rho\left(  w^{\prime}\right)  \right\vert .
\end{align*}
Now for $I\in\mathcal{G}_{\lambda}\left[  S\right]  $ define the right hand
side of (\ref{exclaim}) to be
\begin{align*}
&  w_{I}^{a}\left(  f\right)  \equiv\int_{\mathbb{R}^{3}}\left\vert
T_{I}f\left(  z\right)  \right\vert \zeta_{\lambda}\left(  z-a\right)
dz=\int_{\mathbb{R}^{3}}\left\vert T_{I}f\left(  z\right)  \right\vert
\zeta\left(  \frac{z-a}{2^{\lambda}}\right)  \frac{dz}{2^{3\lambda}}\\
&  =\int_{\mathbb{R}^{3}}\left\vert \widehat{f_{I}^{\Phi}}\left(  z\right)
\right\vert \zeta\left(  \frac{z-a}{2^{\lambda}}\right)  \frac{dz}%
{2^{3\lambda}}\approx\frac{1}{\left\vert B\left(  a,2^{\lambda}\right)
\right\vert }\int_{B\left(  a,2^{\lambda}\right)  }\left\vert \widehat
{f_{I}^{\Phi}}\left(  z\right)  \right\vert \ ,
\end{align*}
and refer to $w_{I}^{a}\left(  f\right)  $ as the `weight' of $\widehat
{f_{I}^{\Phi}}$ relative to the ball $B\left(  a,2^{\lambda}\right)  $, which
represents that portion of the integral of $\left\vert \widehat{f_{I}^{\Phi}%
}\left(  z\right)  \right\vert $ that is taken over the ball $B\left(
a,2^{\lambda}\right)  $. Note that $w_{I}^{a}\left(  f\right)  \lesssim
\left\Vert \widehat{f_{I}^{\Phi}}\right\Vert _{L^{\infty}}\lesssim\left\Vert
\mathbf{1}_{I}f\right\Vert _{L^{1}}\leq\left\vert I\right\vert \left\Vert
f\right\Vert _{L^{\infty}}=2^{-2\lambda}$.

Summarizing, we have%
\begin{equation}
\left\vert T_{I}f\left(  \xi\right)  \right\vert \leq\int_{\mathbb{R}^{3}%
}\left\vert T_{I}f\left(  z\right)  \right\vert \zeta_{\lambda}\left(
z-a\right)  dz=w_{I}^{a}\left(  f\right)  \ ,\ \ \ \ \ \text{for }\xi\in
B\left(  a,2^{\lambda}\right)  .\label{pre exclaim'}%
\end{equation}
and
\begin{equation}
\int_{\mathbb{R}^{3}}\left\vert T_{I}f\left(  z\right)  \right\vert
\zeta_{\lambda}\left(  z-\xi\right)  dz\approx w_{I}^{a}\left(  f\right)
\text{ },\ \ \ \ \ \text{for }\xi\in B\left(  a,2^{\lambda}\right)
,\label{exclaim'}%
\end{equation}
since for $\xi\in B\left(  a,2^{\lambda}\right)  $ we have $\zeta_{\lambda
}\left(  \xi-z\right)  \approx\zeta_{\lambda}\left(  a-z\right)
=\zeta_{\lambda}\left(  z-a\right)  $ by (\ref{pre note that}) and
(\ref{note that}).

Now set%
\[
w_{\ast}^{a}\left(  f\right)  \equiv\max_{I\in\mathcal{G}_{\lambda}\left[
S\right]  }w_{I}^{a}\left(  f\right)  =\max_{I\in\mathcal{G}_{\lambda}\left[
S\right]  }\int_{\mathbb{R}^{3}}\left\vert T_{I}f\left(  z\right)  \right\vert
\zeta_{\lambda}\left(  z-a\right)  dz,
\]
and fix $I_{\ast}^{a}$ such that
\[
w_{I_{\ast}^{a}}^{a}=w_{\ast}^{a}.
\]

Now we replace the function $f$ above with a function in $\mathcal{F}%
_{\kappa,\infty}\left(  s\right)  $ for the remainder of the proof, which with
a small abuse of notation we write as $\mathsf{Q}_{\kappa,s,U}f$ with
$\left\Vert f\right\Vert _{L^{\infty}\left(  U\right)  }\leq1$. For
$\mathsf{Q}_{\kappa,s,U}f\in\mathcal{F}_{\kappa,\infty}\left(  s\right)  $,
$\lambda>1$, and $\alpha\in\mathbb{N}$ chosen appropriately, we will estimate
the contributions to the norm $\left\Vert T\mathsf{Q}_{\kappa,s,U}f\right\Vert
_{L^{q}\left(  \mathbb{R}^{3}\right)  }$ in two exhaustive cases in turn. The
first case will yield the growth factor $2^{\varepsilon s}$, while the next
case will be absorbed. In fact, we show at the end of the proof that we may
take%
\[
\alpha=2\text{\ and }\lambda>\frac{3q}{q-3}.
\]

\end{proof}

\subsection{Case 1: Separated interaction\label{Sub Case 1}}

\begin{proof}
[Proof continued]In \textbf{Case 1} we assume the following property. There
exists a triple of squares $I_{0}^{a},J_{0}^{a},K_{0}^{a}\in\mathcal{G}%
_{\lambda}\left[  U\right]  $ such that%
\begin{align*}
w_{I_{0}^{a}},w_{J_{0}^{a}}  &  >2^{-\alpha\lambda}w_{\ast}^{a}%
\ ,\ \ \ \ \ \text{and }\left\vert \mathbf{c}_{I_{0}^{a}}-\mathbf{c}%
_{J_{0}^{a}}\right\vert ,\left\vert \mathbf{c}_{J_{0}^{a}}-\mathbf{c}%
_{K_{0}^{a}}\right\vert ,\left\vert \mathbf{c}_{K_{0}^{a}}-\mathbf{c}%
_{I_{0}^{a}}\right\vert >2^{10}2^{-\lambda}\ ,\\
&  \fbox{$%
\begin{array}
[c]{ccccc}%
\mathbf{c}_{I_{0}^{a}} &  &  &  & \\
& \cdot &  &  & \\
&  & \mathbf{c}_{K_{0}^{a}} & \leftrightarrows & \mathbf{c}_{K_{0}^{a}}\\
&  &  & \cdot & \\
&  &  &  & \mathbf{c}_{J_{0}^{a}}%
\end{array}
$},
\end{align*}
i.e. there exists a `$2^{10}2^{-\lambda}$-separated' triple $I_{0}^{a}%
,J_{0}^{a},K_{0}^{a}$ of squares of side length $2^{-\lambda}$, such that each
of $I_{0}^{a}$, $J_{0}^{a}\ $and $K_{0}^{a}$ have near maximal weight. In
\textbf{Case 1} we will use the $\nu$-disjoint trilinear estimate
$\mathcal{W}_{\operatorname*{disj}\nu}\left(  \otimes_{3}L^{\infty}\rightarrow
L^{\frac{q}{3}};\varepsilon\right)  $ with $\nu=2^{10}2^{-\lambda}$. For
$\xi\in B\left(  a,2^{\lambda}\right)  $ we throw away the unimodular function
$e^{-i\Phi\left(  c_{I}\right)  \cdot\xi}$, and using (\ref{pre exclaim'}), we
estimate that for $\xi\in B\left(  a,2^{\lambda}\right)  $,%

\begin{align}
\left\vert T\mathsf{Q}_{\kappa,s,U}f\left(  \xi\right)  \right\vert  &
=\left\vert \sum_{L\in\mathcal{G}_{\lambda}\left[  U\right]  }e^{i\phi\left(
\xi,c_{I}\right)  }T_{L}\mathsf{Q}_{\kappa,s,U}f\left(  \xi\right)
\right\vert \leq\sum_{L\in\mathcal{G}_{\lambda}\left[  U\right]  }\left\vert
T_{L}\mathsf{Q}_{\kappa,s,U}f\left(  \xi\right)  \right\vert \label{T weight}%
\\
&  \lesssim\sum_{L\in\mathcal{G}_{\lambda}\left[  U\right]  }w_{L}%
^{a}<2^{2\lambda}w_{\ast}^{a}<2^{\left(  2+\alpha\right)  \lambda}\left(
w_{I_{1}^{a}}^{a}w_{I_{2}^{a}}^{a}w_{I_{3}^{a}}^{a}\right)  ^{\frac{1}{3}%
},\nonumber
\end{align}
since the\emph{\ fixed} triple $\left(  I_{0}^{a},J_{0}^{a},K_{0}^{a}\right)
$ satisfies the near maximal weight condition in \textbf{Case 1}:%
\[
w_{\ast}^{a}<\min\left\{  2^{\alpha\lambda}w_{I_{0}^{a}}^{a},2^{\alpha\lambda
}w_{J_{0}^{a}}^{a},2^{\alpha\lambda}w_{K_{0}^{a}}^{a}\right\}  \leq
2^{\alpha\lambda}\left(  w_{I_{0}^{a}}^{a}\right)  ^{\frac{1}{3}}\left(
w_{J_{0}^{a}}^{a}\right)  ^{\frac{1}{3}}\left(  w_{K_{0}^{a}}^{a}\right)
^{\frac{1}{3}}.
\]
For convenience we suppress dependence on $\kappa$ and $U$ and write%
\[
f^{s}\equiv\mathsf{Q}_{\kappa,s,U}f.
\]
Let $\nu=2^{10}2^{-\lambda}$. Then for $q>3$ and $\xi\in B\left(
a,2^{\lambda}\right)  $, we have from (\ref{T weight}) and (\ref{exclaim'}),
followed by H\"{o}lder's inequality, that%
\begin{align*}
&  \left\vert Tf\left(  \xi\right)  \right\vert ^{q}\lesssim2^{q\left(
2+\alpha\right)  \lambda}\left(  w_{I_{0}^{a}}^{a}w_{J_{0}^{a}}^{a}%
w_{K_{0}^{a}}^{a}\right)  ^{\frac{q}{3}}\\
&  \approx2^{q\left(  2+\alpha\right)  \lambda}\left(  \int_{\mathbb{R}^{3}%
}\left\vert T_{I_{0}^{a}}f^{s}\left(  z_{1}\right)  \right\vert \zeta
_{\lambda}\left(  z_{1}-a\right)  dz_{1}\right)  ^{\frac{q}{3}}\left(
\int_{\mathbb{R}^{3}}\left\vert T_{J_{0}^{a}}f^{s}\left(  z_{2}\right)
\right\vert \zeta_{\lambda}\left(  z_{2}-a\right)  dz_{2}\right)  ^{\frac
{q}{3}}\left(  \int_{\mathbb{R}^{3}}\left\vert T_{K_{0}^{a}}f^{s}\left(
z_{3}\right)  \right\vert \zeta_{\lambda}\left(  z_{3}-a\right)
dz_{3}\right)  ^{\frac{q}{3}}\\
&  \lesssim2^{q\left(  2+\alpha\right)  \lambda}\int_{\mathbb{R}^{3}}%
\int_{\mathbb{R}^{3}}\int_{\mathbb{R}^{3}}\left\vert T_{I_{0}^{a}}f^{s}\left(
\xi-z_{1}\right)  T_{J_{0}^{a}}f^{s}\left(  \xi-z_{2}\right)  T_{K_{0}^{a}%
}f^{s}\left(  \xi-z_{3}\right)  \right\vert ^{\frac{q}{3}}\zeta_{\lambda
}\left(  z_{1}\right)  \zeta_{\lambda}\left(  z_{2}\right)  \zeta_{\lambda
}\left(  z_{3}\right)  dz_{1}dz_{2}dz_{3}\\
&  \lesssim2^{q\left(  2+\alpha\right)  \lambda}\sum_{\left(  I,J,K\right)
\in\mathcal{G}_{\lambda}^{\nu-\mathop{\rm separated}}\left[  U\right]  }%
\int_{\mathbb{R}^{9}}\left\vert T_{I}f^{s}\left(  \xi-z_{1}\right)  T_{J}%
f^{s}\left(  \xi-z_{2}\right)  T_{K}f^{s}\left(  \xi-z_{3}\right)  \right\vert
^{\frac{q}{3}}d\mu_{\lambda}\left(  z_{1},z_{2},z_{3}\right)  \ ,
\end{align*}
where
\[
\mathcal{G}_{\lambda}^{\nu-\mathop{\rm separated}}\left[  U\right]
\equiv\left\{  \left(  I,J,K\right)  \in\mathcal{G}_{\lambda}:\left(
I,J,K\right)  \text{ is }\nu\text{-separated as in (\ref{weak sep})}\right\}
,
\]
and%
\[
d\mu_{\lambda}\left(  z_{1},z_{2},z_{3}\right)  \equiv\zeta_{\lambda}\left(
z_{1}\right)  \zeta_{\lambda}\left(  z_{2}\right)  \zeta_{\lambda}\left(
z_{3}\right)  dz_{1}dz_{2}dz_{3}%
\]
is a bounded multiple of a probability measure.

Then we have%
\begin{align*}
&  \int_{B\left(  a,2^{\lambda}\right)  }\left\vert Tf^{s}\left(  \xi\right)
\right\vert ^{q}d\xi\\
&  \lesssim2^{q\left(  2+\alpha\right)  \lambda}\int_{B\left(  a,2^{\lambda
}\right)  }\sum_{\left(  I,J,K\right)  \in\mathcal{G}_{\lambda}^{\nu
-\mathop{\rm separated}}\left[  U\right]  }\int_{\mathbb{R}^{9}}\left\vert
T_{I}f^{s}\left(  \xi-z_{1}\right)  T_{J}f^{s}\left(  \xi-z_{2}\right)
T_{K}f^{s}\left(  \xi-z_{3}\right)  \right\vert ^{\frac{q}{3}}d\mu_{\lambda
}\left(  z_{1},z_{2},z_{3}\right)  d\xi\\
&  =2^{q\left(  2+\alpha\right)  \lambda}\sum_{\left(  I,J,K\right)
\in\mathcal{G}_{\lambda}^{\nu-\mathop{\rm separated}}\left[  U\right]  }%
\int_{\mathbb{R}^{9}}\left\{  \int_{B\left(  a,2^{\lambda}\right)  }\left\vert
T_{I}f^{s}\left(  \xi-z_{1}\right)  T_{J}f^{s}\left(  \xi-z_{2}\right)
T_{K}f^{s}\left(  \xi-z_{3}\right)  \right\vert ^{\frac{q}{3}}d\xi\right\}
d\mu_{\lambda}\left(  z_{1},z_{2},z_{3}\right)  .
\end{align*}

Now consider those $a\in\Gamma_{\lambda}\left(  R\right)  $ for which
\textbf{Case 1} is in effect for the ball $B\left(  a,2^{\lambda}\right)  $
and denote by $\Gamma_{\lambda}\left(  \text{\textbf{Case 1}}\right)  $ the
union of all the balls $B\left(  a,2^{\lambda}\right)  $ for which $a$ is in
\textbf{Case 1}. Summing over points $a\in\Gamma_{\lambda}\left(  R\right)  $
such that \textbf{Case 1} is in effect for the ball $B\left(  a,2^{\lambda
}\right)  $, we obtain%
\begin{align*}
&  \sum_{a\in\Gamma_{\lambda}\left(  \text{\textbf{Case 1}}\right)  }%
\int_{B\left(  a,2^{\lambda}\right)  }\left\vert Tf^{s}\left(  \xi\right)
\right\vert ^{q}d\xi\\
&  \lesssim2^{q\left(  2+\alpha\right)  \lambda}\sum_{\left(  I,J,K\right)
\in\mathcal{G}_{\lambda}^{\nu-\mathop{\rm separated}}\left[  U\right]  }%
\int_{\mathbb{R}^{9}}\left\{  \sum_{a\in\Gamma_{\lambda}\left(
\text{\textbf{Case 1}}\right)  }\int_{B\left(  a,2^{\lambda}\right)
}\left\vert T_{I}f^{s}\left(  \xi-z_{1}\right)  T_{J}f^{s}\left(  \xi
-z_{2}\right)  T_{K}f^{s}\left(  \xi-z_{3}\right)  \right\vert ^{\frac{q}{3}%
}d\xi\right\}  d\mu_{\lambda}\left(  z_{1},z_{2},z_{3}\right) \\
&  \lesssim2^{q\left(  2+\alpha\right)  \lambda}\sum_{\left(  I,J,K\right)
\in\mathcal{G}_{\lambda}^{\nu-\mathop{\rm separated}}\left[  U\right]  }%
\int_{\mathbb{R}^{9}}\left\{  \int_{B_{R}}\left\vert T_{I}f^{s}\left(
\xi-z_{1}\right)  T_{J}f^{s}\left(  \xi-z_{2}\right)  T_{K}f^{s}\left(
\xi-z_{3}\right)  \right\vert ^{\frac{q}{3}}d\xi\right\}  d\mu_{\lambda
}\left(  z_{1},z_{2},z_{3}\right)  .
\end{align*}
This is now dominated by%
\begin{align*}
&  2^{q\left(  2+\alpha\right)  \lambda}\left(  \#\mathcal{G}_{\lambda}%
^{\nu-\mathop{\rm separated}}\left[  U\right]  \right)  ^{3}\int
_{\mathbb{R}^{9}}\left\{  C_{\varepsilon}^{q}2^{q\varepsilon s}\left\Vert
f_{1}\right\Vert _{L^{\infty}}\left\Vert f_{2}\right\Vert _{L^{\infty}%
}\left\Vert f_{3}\right\Vert _{L^{\infty}}\right\}  d\mu_{\lambda}\left(
z_{1},z_{2}...,z_{N}\right) \\
&  \leq2^{q\left(  2+\alpha\right)  \lambda}\left(  \#\mathcal{G}_{\lambda
}^{\nu-\mathop{\rm separated}}\left[  U\right]  \right)  ^{3}\int
_{\mathbb{R}^{9}}\left\{  C_{\varepsilon}^{q}2^{q\varepsilon s}\left(
2^{\frac{2}{r}s}\left\Vert f\right\Vert _{L^{\infty}}\right)  ^{3}\right\}
d\mu_{\lambda}\left(  z_{1},z_{2}...,z_{N}\right)  \lesssim C_{\varepsilon
}^{q}2^{q\left(  2+\alpha\right)  \lambda}2^{6\lambda}2^{q\varepsilon
s}2^{\frac{6}{r}s},
\end{align*}

by (\ref{penalty}), which gives $\left\Vert f_{k}\right\Vert _{L^{\infty}}\leq
C_{r}2^{\frac{2}{r}s}\left\Vert f\right\Vert _{L^{\infty}}=C_{r}2^{\frac{2}%
{r}s}$ here, upon appealing to the local single scale $\nu$-disjoint trilinear
assumption $\mathcal{A}_{\operatorname*{disj}\nu}^{s,\kappa,\delta}\left(
\otimes_{3}L^{\infty}\rightarrow L^{\frac{q}{3}};\varepsilon\right)  $ with
$\nu=2^{10}2^{-\lambda}$, and with%
\[
f_{1}=\widetilde{\mathsf{M}_{z_{1}}}f_{I}^{s},\ \ \ f_{2}=\widetilde
{\mathsf{M}_{z_{1}}}f_{J}^{s},\ \ \ \ f_{3}=\widetilde{\mathsf{M}_{z_{1}}%
}f_{K}^{s},
\]
where $\widetilde{\mathsf{M}_{z}}\left(  x\right)  =e^{i\left\langle
z,\Phi\left(  x\right)  \right\rangle }$. Indeed, recall from \cite{RiSa} that
if $\mathsf{M}_{z}\left(  w\right)  =e^{i\left\langle z,w\right\rangle }$,
then $\mathsf{M}_{z}\Phi_{\ast}=\Phi_{\ast}\widetilde{\mathsf{M}_{z}}$ and
thus,%
\begin{align*}
\left\vert T_{I}f^{s}\left(  \xi-z_{1}\right)  \right\vert  &  =\widehat
{\Phi_{\ast}f_{I}^{s}}\left(  \xi-z_{1}\right)  =\widehat{\mathsf{M}_{z_{1}%
}\Phi_{\ast}f_{I}^{s}}\left(  \xi\right)  =\widehat{\Phi_{\ast}\widetilde
{\mathsf{M}_{z_{1}}}f_{I}^{s}}\left(  \xi\right)  =\mathcal{E}f_{1}\left(
\xi\right)  ,\\
\text{and }\left\vert T_{J}f^{s}\left(  \xi-z_{2}\right)  \right\vert  &
=\mathcal{E}f_{1}\left(  \xi\right)  \text{ and }\left\vert T_{K}f^{s}\left(
\xi-z_{1}\right)  \right\vert =\mathcal{E}f_{3}\left(  \xi\right)  \text{,}%
\end{align*}
and%
\begin{align*}
&  \int_{B_{\delta}\left(  0,2^{s}\right)  }\left\vert T_{I}f^{s}\left(
\xi-z_{1}\right)  T_{J}f^{s}\left(  \xi-z_{2}\right)  T_{K}f^{s}\left(
\xi-z_{3}\right)  \right\vert ^{\frac{q}{3}}d\xi=\int_{B_{\delta}\left(
0,2^{s}\right)  }\left\vert \mathcal{E}f_{1}\left(  \xi\right)  \mathcal{E}%
f_{2}\left(  \xi\right)  \mathcal{E}f_{3}\left(  \xi\right)  \right\vert
^{\frac{q}{3}}d\xi\\
&  =\left\Vert \prod_{j=1}^{3}\mathcal{E}_{j}f_{j}\right\Vert _{L^{^{\frac
{q}{3}}}\left(  B_{\delta}\left(  0,2^{s}\right)  \right)  }^{\frac{q}{3}}%
\leq\left(  C_{\varepsilon,\nu,q}R^{\varepsilon}\right)  ^{q}\prod_{j=1}%
^{3}\left\Vert f_{j}\right\Vert _{L^{\infty}}^{\frac{q}{3}}\lesssim\left(
C_{\varepsilon,\nu,q,r}2^{\varepsilon s}\right)  ^{q}2^{\frac{2q}{r}s},
\end{align*}
since the triple $\left(  I,J,K\right)  $ is $\nu$-separated as in
(\ref{weak sep}), and since $\left\vert f_{j}\right\vert \leq1$. As a
consequence we have%
\begin{align*}
&  \sum_{\left(  I,J,K\right)  \in\mathcal{G}_{\lambda}^{\nu
-\mathop{\rm separated}}\left[  U\right]  }\int_{B_{\delta}\left(
0,2^{s}\right)  }\left\vert T_{I}f^{s}\left(  \xi\right)  T_{J}f^{s}\left(
\xi\right)  T_{K}f^{s}\left(  \xi\right)  \right\vert ^{\frac{q}{3}}d\xi\\
&  \leq\sum_{\left(  I,J,K\right)  \in\mathcal{G}_{\lambda}^{\nu
-\mathop{\rm separated}}\left[  U\right]  }\left(  C_{\varepsilon,\nu
,q,r}2^{\varepsilon s}\right)  ^{q}2^{\frac{2q}{r}s}\lesssim2^{6\lambda
}\left(  C_{\varepsilon,\nu,q}2^{\varepsilon s}\right)  ^{q}2^{\frac{2q}{r}%
s}=\left(  C_{\varepsilon,\nu,q,r}\right)  ^{q}2^{6\lambda}2^{q\varepsilon
s}2^{\frac{2q}{r}s}.
\end{align*}
Altogether then, we have proved that
\[
\left\Vert \mathbf{1}_{\Gamma_{\lambda}\left(  \text{\textbf{Case 1}}\right)
}T\mathsf{Q}_{\kappa,s,U}f\right\Vert _{L^{q}\left(  B_{\delta}\left(
0,2^{s}\right)  \right)  }=\left\Vert \mathbf{1}_{\Gamma_{\lambda}\left(
\text{\textbf{Case 1}}\right)  }Tf^{s}\right\Vert _{L^{q}\left(  B_{\delta
}\left(  0,2^{s}\right)  \right)  }\lesssim C_{\varepsilon,\nu,q,r}2^{\left(
\frac{6}{q}+2+\alpha\right)  \lambda}2^{\left(  \varepsilon+\frac{2}%
{r}\right)  s},
\]
where $\mathbf{1}_{\Gamma_{\lambda}\left(  \text{\textbf{Case 1}}\right)  }$
indicates the union of those balls $B\left(  a,2^{\lambda}\right)  $ for which
\textbf{Case 1} holds. We can absorb the small constant $\frac{2}{r}$ into
$\varepsilon$ simply by taking $r$ large in (\ref{penalty}), and the bound we
will carry forward to the end of the proof will be $C_{\varepsilon,\nu
,q}2^{\left(  \frac{6}{q}+2+\alpha\right)  \lambda}2^{\varepsilon s}$.
\end{proof}

\begin{remark}
A key idea used in \textbf{Case 1} is the domination of $\left\vert
T_{I}f\left(  \xi\right)  \right\vert $ by $w_{I}^{a}\left(  f\right)  \ $for
$\xi\in B\left(  a,2^{\lambda}\right)  $ in (\ref{pre exclaim'}) in order to
permit the pigeonholing. This domination must\ then be `undone' by
H\"{o}lder's inequality with exponent $\frac{q}{3}$, in order to appeal to the
$\nu$-disjoint trilinear assumption $\mathcal{E}_{\mathop{\rm disj}\nu}\left(
\otimes_{3}L^{\infty}\rightarrow L^{\frac{q}{3}}\right)  $. This feature of
the argument is an obstruction to generalizations to $N$-linear extensions in
higher dimensions $n$, since the critical exponent for $q$, namely $\frac
{2n}{n-1}$, is strictly less than $3\leq N$ when $n\geq4$, and so we cannot
assume $\frac{q}{N}\geq1$.
\end{remark}

\subsection{Case 2: Clustered interaction\label{Sub Case 2}}

\begin{proof}
[Proof continued]In \textbf{Case 2} we assume the negation of \textbf{Case 1}.
Thus if $\left\{  I_{1}^{a},I_{2}^{a}\right\}  $ is a maximal pair of cubes
such that $w_{I_{1}^{a}},w_{I_{2}^{a}}>2^{-\alpha\lambda}w_{\ast}^{a}\ $and
$\left\vert \mathbf{c}_{I_{1}^{a}}-\mathbf{c}_{I_{2}^{a}}\right\vert
>2^{10}2^{-\lambda}$, then
\begin{equation}
w_{I}\leq2^{-\alpha\lambda}w_{\ast}^{a},\ \ \ \ \ \text{ if }%
\operatorname*{dist}\left(  I,I_{1}^{a}\cup I_{2}^{a}\right)  >2^{10}%
2^{-\lambda}.\label{small away}%
\end{equation}
In this case we will use rescaling and recursion as in \cite{TaVaVe}.

Let $\xi\in B\left(  a,2^{\lambda}\right)  $ for some $a\in\Gamma_{\lambda
}\left(  R\right)  $. Using (\ref{small away}), and with $\left\vert
\cdot\right\vert _{\operatorname*{square}}$ denoting the `square' norm in
$\mathbb{R}^{3}$,
\begin{align*}
\left\vert Tf^{s}\left(  \xi\right)  \right\vert  &  =\left\vert \sum
_{I\in\mathcal{G}_{\lambda}\left[  U\right]  }\int_{I}e^{i\phi\left(
\xi,y\right)  }f^{s}\left(  y\right)  dy\right\vert \\
&  \leq\sum_{m=1}^{2}\left\vert \sum_{I\in\mathcal{G}_{\lambda}\left[
U\right]  :\ \left\vert c_{I}-c_{I_{m}^{a}}\right\vert
_{\operatorname*{square}}\leq2^{10}2^{-\lambda}}\int_{I}e^{i\phi\left(
\xi,y\right)  }f^{s}\left(  y\right)  dy\right\vert +\sum_{I\in\mathcal{G}%
_{\lambda}\left[  U\right]  :\ \operatorname*{dist}\left(  I,\bigcup_{m=1}%
^{M}I_{m}^{a}\right)  >2^{10}2^{-\lambda}}\left\vert T_{I}f^{s}\left(
\xi\right)  \right\vert \\
&  \leq C\max_{K\in\mathcal{G}_{\lambda-10}\left[  U\right]  }\left\vert
\int_{K}e^{i\phi\left(  \xi,y\right)  }f^{s}\left(  y\right)  dy\right\vert
+\sum_{I\in\mathcal{G}_{\lambda}\left[  U\right]  :\ \operatorname*{dist}%
\left(  I,\bigcup_{m=1}^{2}I_{m}^{a}\right)  >2^{10}2^{-\lambda}}\left\vert
T_{I}f^{s}\left(  \xi\right)  \right\vert \\
&  \leq C\max_{K\in\mathcal{G}_{\lambda-10}\left[  U\right]  }\left\vert
T_{K}f^{s}\left(  \xi\right)  \right\vert +\left(  \#\mathcal{G}_{\lambda
}\left[  U\right]  \right)  2^{-\alpha\lambda}w_{\ast}^{a}\leq C\max
_{K\in\mathcal{G}_{\lambda-10}\left[  U\right]  }\left\vert T_{K}f^{s}\left(
\xi\right)  \right\vert +2^{\left(  2-\alpha\right)  \lambda}w_{\ast}^{a}\ ,
\end{align*}
since if $K_{m}^{a}\in\mathcal{G}_{\lambda-10}\left[  U\right]  $ is the
unique such cube containing $I_{m}^{a}$, then we decompose $f^{s}=f_{K_{m}%
^{a}}^{s}+\sum_{I\in\mathcal{G}_{\lambda}\left[  S\right]  :\ I\cap K_{m}%
^{a}=\emptyset}f_{I}^{s}$, and without loss of generality we may also assume
$\left\vert c_{I}-c_{I_{m}^{a}}\right\vert \geq2^{10-\lambda}$. Now
\begin{align*}
\int\left\vert T_{I}f^{s}\left(  \xi-z\right)  \right\vert \zeta_{\lambda}%
^{a}\left(  z\right)  dz  &  \leq\left(  \int\left\vert T_{I}f^{s}\left(
\xi-z\right)  \right\vert ^{q}\zeta_{\lambda}\left(  z-a\right)  dz\right)
^{\frac{1}{q}}\left(  \int\zeta_{\lambda}\left(  z-a\right)  dz\right)
^{\frac{1}{q^{\prime}}}\\
&  \lesssim\left(  \int\left\vert T_{I}f^{s}\left(  z-a\right)  \right\vert
^{q}\zeta_{\lambda}\left(  z\right)  dz\right)  ^{\frac{1}{q}},
\end{align*}
and so for $\xi\in B\left(  a,2^{\lambda}\right)  $,
\begin{align*}
\left\vert Tf^{s}\left(  \xi\right)  \right\vert ^{q}  &  \leq C\sum
_{K\in\mathcal{G}_{\lambda-10}\left[  U\right]  }\left\vert T_{K}f^{s}\left(
\xi\right)  \right\vert ^{q}+C2^{\left(  2-\alpha\right)  \lambda q}%
\int\left\vert T_{I_{\ast}^{a}}f^{s}\left(  z\right)  \right\vert ^{q}%
\zeta_{\lambda}\left(  z-a\right)  dz\\
&  \leq C\sum_{K\in\mathcal{G}_{\lambda-10}\left[  U\right]  }\left\vert
T_{K}f^{s}\left(  \xi\right)  \right\vert ^{q}+C2^{\left(  2-\alpha\right)
\lambda q}\sum_{I\in\mathcal{G}_{\lambda}\left[  U\right]  }\int\left\vert
T_{I}f^{s}\left(  z\right)  \right\vert ^{q}\zeta_{\lambda}\left(  z-a\right)
dz,
\end{align*}
where $C$ is independent of all the relevant parameters, and where we have
added back in all $K\in\mathcal{G}_{\lambda^{\prime}}\left[  S\right]  $
rather than just $K_{\ast}^{a}$, and all $I\in\mathcal{G}_{\lambda}\left[
S\right]  $ rather than just $I_{\ast}^{a}$, in order that the dominating
expression on the right hand side is indepenent of $a\in\Gamma_{\lambda
}\left(  R\right)  $.

Summing over $a\in\Gamma_{\lambda}\left(  R\right)  $, we see that the
corresponding contribution over $B_{R}$ is at most%
\begin{align}
& \label{contrib Case 2}\\
\left\Vert \mathbf{1}_{\Gamma_{\lambda}\left(  \text{\textbf{Case 2}}\right)
}Tf^{s}\right\Vert _{L^{q}\left(  B_{\delta}\left(  0,2^{s}\right)  \right)
}^{q}  &  \equiv\sum_{a\in\Gamma_{\lambda}}\left\{  C\sum_{K\in\mathcal{G}%
_{\lambda-10}\left[  U\right]  }\int_{B\left(  a,2^{\lambda}\right)
}\left\vert T_{K}f^{s}\left(  \xi\right)  \right\vert ^{q}d\xi+C2^{\left(
2-\alpha\right)  \lambda q}\sum_{I\in\mathcal{G}_{\lambda}\left[  U\right]
}\int\left\vert T_{I}f^{s}\left(  z\right)  \right\vert ^{q}\zeta_{s}%
^{a}\left(  z\right)  dz\right\} \nonumber\\
&  \lesssim C\sum_{K\in\mathcal{G}_{\lambda-10}\left[  U\right]  }%
\int_{\mathbb{R}^{3}}\left\vert T_{K}f^{s}\left(  \xi\right)  \right\vert
^{q}d\xi+C2^{-n\lambda}2^{-\left(  2-\alpha\right)  \lambda q}\sum
_{I\in\mathcal{G}_{\lambda}\left[  U\right]  }\int_{\mathbb{R}^{3}}\left\vert
T_{I}f^{s}\left(  \xi\right)  \right\vert ^{q}d\xi,\nonumber
\end{align}
since $\sum_{a\in\Gamma_{\lambda}}\zeta_{\lambda}\left(  z-a\right)
\lesssim2^{-3\lambda}\mathbf{1}_{\mathbb{R}^{3}}\left(  z\right)
+\operatorname*{rapid}\operatorname*{decay}$.

At this point we follow \cite{BoGu} in using parabolic rescaling, as
introduced in Tao, Vargas and Vega \cite{TaVaVe}, on the integral
\[
\operatorname*{Int}_{\rho}\left(  \xi\right)  \equiv\int_{\left\vert
y-\overline{y}\right\vert <\rho}e^{i\phi\left(  \xi,y\right)  }f^{s}\left(
y\right)  dy=\int_{\left\vert y-\overline{y}\right\vert <\rho}e^{i\left[
\xi_{1}y_{1}+\xi_{2}y_{2}+\xi_{3}\left(  y_{1}^{2}+y_{2}^{2}\right)  \right]
}f^{s}\left(  y\right)  dy,\ \ \ \ \ \text{for }0<\rho<1,
\]
to obtain%
\begin{align*}
&  \left\vert \operatorname*{Int}_{\rho}\left(  \xi\right)  \right\vert
\overset{y=\overline{y}+y^{\prime}}{=}\left\vert \int_{\left\vert y^{\prime
}\right\vert <\rho}e^{i\left[  \xi_{1}\left(  \overline{y}_{1}+y_{1}^{\prime
}\right)  +\xi_{2}\left(  \overline{y}_{2}+y_{2}^{\prime}\right)  +\xi
_{3}\left(  \left(  \overline{y}_{1}+y_{1}^{\prime}\right)  ^{2}+\left(
\overline{y}_{2}+y_{2}^{\prime}\right)  ^{2}\right)  \right]  }f^{s}\left(
\overline{y}+y^{\prime}\right)  dy^{\prime}\right\vert \\
&  =\left\vert \int_{\left\vert y^{\prime}\right\vert <\rho}e^{i\left[
\left(  \xi_{1}+2\overline{y_{1}}\xi_{3}\right)  y_{1}^{\prime}+\left(
\xi_{2}+2\overline{y_{2}}\xi_{3}\right)  y_{2}^{\prime}+\xi_{3}\left\vert
y^{\prime}\right\vert ^{2}\right]  }f^{s}\left(  \overline{y}+y^{\prime
}\right)  dy^{\prime}\right\vert .
\end{align*}
Thus we conclude that%
\begin{align}
&  \ \ \ \ \ \ \ \ \ \ \ \ \ \ \ \ \ \ \ \ \ \ \ \ \ \ \ \ \ \ \left\Vert
\operatorname*{Int}_{\rho}\right\Vert _{L^{q}\left(  B_{R}\right)  }=\left(
\int_{B_{R}}\left\vert \operatorname*{Int}_{\rho}\left(  \xi\right)
\right\vert ^{q}d\xi\right)  ^{\frac{1}{q}}\label{Int est}\\
&  =\left(  \int_{B_{R}}\left\vert \int_{\left\vert y^{\prime}\right\vert
<\rho}e^{i\left[  \left(  \xi_{1}+2\overline{y_{1}}\xi_{3}\right)
y_{1}^{\prime}+...+\left(  \xi_{2}+2\overline{y_{2}}\xi_{3}\right)
y_{2}^{\prime}+\xi_{3}\left\vert y^{\prime}\right\vert ^{2}\right]  }%
f^{s}\left(  \overline{y}+y^{\prime}\right)  dy^{\prime}\right\vert ^{q}%
d\xi\right)  ^{\frac{1}{q}}\nonumber\\
&  =\left(  \int_{B_{R}}\left\vert \int_{\left\vert y^{\prime}\right\vert
<\rho}e^{i\left[  \left(  \rho\xi_{1}+2\overline{\frac{y_{1}}{\rho}}\rho
^{2}\xi_{3}\right)  \frac{y_{1}^{\prime}}{\rho}+\left(  \rho\xi_{2}%
+2\overline{\frac{y_{2}}{\rho}}\rho^{2}\xi_{3}\right)  y_{2}^{\prime}+\rho
^{2}\xi_{3}\left\vert \frac{y^{\prime}}{\rho}\right\vert ^{2}\right]  }%
f^{s}\left(  \overline{y}+y^{\prime}\right)  \rho^{2}d\left(  \frac{y^{\prime
}}{\rho}\right)  \right\vert ^{q}\frac{d\left(  \rho\xi^{\prime}\right)
d\left(  \rho^{2}\xi_{3}\right)  }{\rho^{4}}\right)  ^{\frac{1}{q}}\nonumber\\
&  =\rho^{2}\rho^{-\frac{4}{q}}\left(  \int_{B_{\rho R}}\left\vert
\int_{\left\vert y^{\prime}\right\vert <1}e^{i\left[  \left(  \xi
_{1}+2\overline{y_{1}}\xi_{3}\right)  y_{1}^{\prime}+\left(  \xi
_{2}+2\overline{y_{2}}\xi_{3}\right)  y_{2}^{\prime}+\xi_{3}\left\vert
y^{\prime}\right\vert ^{2}\right]  }f^{s}\left(  \rho\left(  \overline
{y}+y^{\prime}\right)  \right)  dy^{\prime}\right\vert ^{q}d\xi^{\prime}%
d\xi_{3}\right)  ^{\frac{1}{q}}\leq C\rho^{2}\rho^{-\frac{4}{q}}%
P_{\kappa,\delta^{\prime}}^{\left(  q\right)  }\left(  s-m\right)  ,\nonumber
\end{align}
by using that $2^{\frac{s-m}{1-\delta^{\prime}}}=2^{-m}2^{\frac{s}{1-\delta}}
$ implies
\[
\rho=2^{-m}\text{ and }\delta^{\prime}=\frac{\delta}{1-\frac{m}{s}\left(
1-\delta\right)  }\text{ and }s>m,
\]
and using Corollary \ref{F inv} applied to%
\[
f^{s}\left(  \rho\left(  \overline{y}+y^{\prime}\right)  \right)
=\mathsf{Q}_{\kappa,s,U}f\left(  2^{-m}\left(  \overline{y}+y^{\prime}\right)
\right)
\]
which says that $f^{s}\left(  \rho\left(  \overline{y}+y^{\prime}\right)
\right)  $ belongs to $\mathcal{F}_{\kappa,\infty}\left(  s-m\right)  $ for
$\rho=2^{-m}$, $m\geq0$. Note that we have used%
\[
\rho R=2^{-m}2^{\frac{s}{1-\delta}}=2^{\frac{s-m}{1-\delta^{\prime}}}\text{
implies }\delta^{\prime}=\frac{\delta}{1-\frac{m}{s}\left(  1-\delta\right)
},
\]

and that if $s>m$, then $\delta<\delta^{\prime}<1$.

We also used that
\[
\delta_{m,\infty}f^{s}\left(  x\right)  \equiv f^{s}\left(  2^{-m}x\right)
=f^{s}\left(  \rho x\right)
\]
for $\rho=2^{-m}$ has the same $L^{\infty}$ norm as $f^{s}$, and finally that
the paraboloid is invariant under parabolic rescaling. Note that the factor
$\rho^{2}$ arises from $\left\vert y^{\prime}\right\vert <\rho$, and that the
factor $\rho^{-\frac{4}{q}}$ arises from parabolic rescaling. These features
remain in play for an arbitrary quadratic surface of positive Gaussian curvature.

Thus using (\ref{Int est}), first with $\rho=2^{10-\lambda}$ and then with
$\rho=2^{-\lambda}$, together with (\ref{contrib Case 2}), we obtain that the
contribution $\left\Vert \mathbf{1}_{\Gamma_{\lambda}\left(
\text{\textbf{Case 2}}\right)  }Tf^{s}\right\Vert _{L^{q}\left(  B_{\delta
}\left(  0,2^{s}\right)  \right)  }^{q}$ to the norm $\left\Vert
Tf^{s}\right\Vert _{L^{q}\left(  B_{\delta}\left(  0,2^{s}\right)  \right)  }$
satisfies:
\begin{align}
&  \left\Vert \mathbf{1}_{\Gamma_{\lambda}\left(  \text{\textbf{Case 2}%
}\right)  }T\mathsf{Q}_{\kappa,s,U}f\right\Vert _{L^{q}\left(  B_{\delta
}\left(  0,2^{s}\right)  \right)  }=\left\Vert \mathbf{1}_{\Gamma_{\lambda
}\left(  \text{\textbf{Case 2}}\right)  }Tf^{s}\right\Vert _{L^{q}\left(
B_{\delta}\left(  0,2^{s}\right)  \right)  }\label{case 2 est}\\
&  \leq C\left(  \#\mathcal{G}_{\lambda-10}\left[  U\right]  \right)
^{\frac{1}{q}}\left(  2^{10-\lambda}\right)  ^{2-\frac{4}{q}}P_{\kappa
,\delta^{\prime}}^{\left(  q\right)  }\left(  s+10-\lambda\right) \nonumber\\
&  +C\left(  \#\mathcal{G}_{\lambda}\left[  U\right]  \right)  ^{\frac{1}{q}%
}2^{-\left(  2-\alpha+\frac{3}{q}\right)  \lambda}\left(  2^{-\lambda}\right)
^{2-\frac{4}{q}}P_{\kappa,\delta^{\prime}}^{\left(  q\right)  }\left(
s-\lambda\right) \nonumber\\
&  =C2^{\left(  \frac{6}{q}-2\right)  \lambda}P_{\kappa,\delta^{\prime}%
}^{\left(  q\right)  }\left(  s+10-\lambda\right)  +C2^{\left(  \frac{3}%
{q}+\alpha-4\right)  \lambda}P_{\kappa,\delta^{\prime}}^{\left(  q\right)
}\left(  s-\lambda\right)  \ ,\nonumber
\end{align}
for $s>\lambda$, since%
\[
\left(  \#\mathcal{G}_{\lambda-10}\left[  S\right]  \right)  ^{\frac{1}{q}%
}\left(  2^{-\left(  \lambda-10\right)  }\right)  ^{2-\frac{4}{q}}=2^{\frac
{2}{q}\left(  \lambda-10\right)  }2^{-2\left(  \lambda-10\right)  }2^{\frac
{4}{q}\left(  \lambda-10\right)  }=2^{\left(  \frac{6}{q}-2\right)  \left(
\lambda-10\right)  },
\]
and%
\[
2^{-\frac{3}{q}\lambda}2^{-\left(  2-\alpha\right)  \lambda}\left(
\#\mathcal{G}_{\lambda}\left[  U\right]  \right)  ^{\frac{1}{q}}\left(
2^{-\lambda}\right)  ^{2-\frac{4}{q}}=2^{-\frac{3}{q}\lambda}2^{-\left(
2-\alpha\right)  \lambda}2^{\frac{2}{q}\lambda}2^{-2\lambda}2^{\frac{4}%
{q}\lambda}=2^{\left(  \frac{3}{q}+\alpha-4\right)  \lambda}\ .
\]

\end{proof}

\subsection{Completing the proof of Theorem \ref{Loc lin}}

\begin{proof}
[Proof continued]Recall from (\ref{Q_R}) that
\[
P_{\kappa,\delta}^{\left(  q\right)  }\left(  s\right)  \equiv\sup
_{\mathcal{G}\in\mathbb{G}}\sup_{\left\Vert f\right\Vert _{L^{\infty}\left(
U\right)  }\leq1}\left(  \int_{B_{\delta}\left(  0,2^{s}\right)  }\left\vert
\mathcal{E}\left(  \Phi_{\ast}\mathsf{Q}_{\kappa,s,U}^{\eta,\mathcal{G}%
}f\right)  \left(  \xi\right)  \right\vert ^{q}d\xi\right)  ^{\frac{1}{q}}.
\]
We now define the quantity%
\[
P_{\kappa,\delta,\ast}^{\left(  q\right)  }\left(  s\right)  \equiv
\sup_{\mathcal{G}\in\mathbb{G}}\sup_{\left\Vert f\right\Vert _{L^{\infty
}\left(  U\right)  }\leq1}\left(  \int_{\mathbb{R}^{3}\setminus B_{\delta
}\left(  0,2^{s}\right)  }\left\vert \mathcal{E}\left(  \Phi_{\ast}%
\mathsf{Q}_{\kappa,s,U}^{\eta,\mathcal{G}}f\right)  \left(  \xi\right)
\right\vert ^{q}d\xi\right)  ^{\frac{1}{q}},
\]
where we are integrating \textbf{away} from the ball $B\left(  0,2^{s}\right)
$. As a consequence, we can use the estimates in \cite[see Case 2 of the proof
of Theorem 10]{RiSa} to show that%
\[
P_{\kappa,\delta,\ast}^{\left(  q\right)  }\left(  s\right)  \lesssim
C_{\kappa,q,p,\delta,N}2^{-N\delta s}\leq C_{\kappa,q,p,\delta,N},
\]
for some large integer $N$ (arising from integration by parts in \cite{RiSa}).
So far, remembering that we have absorbed the penalty factor $2^{\frac{2}{r}%
s}$ into $2^{\varepsilon s}$ in the \textbf{Case 1} bound, we have shown that
for $s>\lambda$, we have%
\begin{align*}
&  \left\Vert T\mathsf{Q}_{\kappa,s,U}f\right\Vert _{L^{q}\left(  B_{\delta
}\left(  0,2^{s}\right)  \right)  }\leq\left\Vert \mathbf{1}_{\Gamma_{\lambda
}\left(  \text{\textbf{Case 1}}\right)  }T\mathsf{Q}_{\kappa,s,U}f\right\Vert
_{L^{q}\left(  B_{\delta}\left(  0,2^{s}\right)  \right)  }+\left\Vert
\mathbf{1}_{\Gamma_{\lambda}\left(  \text{\textbf{Case 2}}\right)
}T\mathsf{Q}_{\kappa,s,U}f\right\Vert _{L^{q}\left(  B_{\delta}\left(
0,2^{s}\right)  \right)  }\\
&  \lesssim C_{\varepsilon}2^{\left(  2+\alpha\right)  \lambda}2^{\frac{6}%
{q}\lambda}2^{\varepsilon s}+C2^{\left(  \frac{6}{q}-2\right)  \lambda
}P_{\kappa,2\delta}^{\left(  q\right)  }\left(  s+10-\lambda\right)
+C2^{\left(  \frac{3}{q}+\alpha-4\right)  \lambda}P_{\kappa,\delta^{\prime}%
}^{\left(  q\right)  }\left(  s-\lambda\right) \\
&  \lesssim C_{\varepsilon,N}2^{\left(  2+\alpha\right)  \lambda}2^{\frac
{6}{q}\lambda}2^{\varepsilon s}+C\left[  2^{\left(  \frac{6}{q}-2\right)
\lambda}+2^{\left(  \frac{3}{q}+\alpha-4\right)  \lambda}\right]  \max_{1\leq
s^{\prime}\leq s+10-\lambda}P_{\kappa,\delta^{\prime}}^{\left(  q\right)
}\left(  s^{\prime}\right) \\
&  \lesssim C_{\varepsilon}2^{\left(  2+\alpha\right)  \lambda}2^{\frac{6}%
{q}\lambda}2^{\varepsilon s}+\frac{1}{2}\max_{1\leq s^{\prime}\leq\lambda
}P_{\kappa,\delta}^{\left(  q\right)  }\left(  s^{\prime}\right)  +\frac{1}%
{2}\max_{\lambda<s^{\prime}\leq s+10-\lambda}P_{\kappa,\delta^{\prime}%
}^{\left(  q\right)  }\left(  s^{\prime}\right) \\
&  \lesssim C_{\varepsilon}2^{\left(  2+\alpha\right)  \lambda}2^{\frac{6}%
{q}\lambda}2^{\varepsilon s}+C_{\lambda,\kappa,\delta,q}+\frac{1}{2}%
\max_{\lambda\leq s^{\prime}\leq s+10-\lambda}P_{\kappa,\delta}^{\left(
q\right)  }\left(  s^{\prime}\right)  +\frac{1}{2}\max_{\lambda<s^{\prime}\leq
s+10-\lambda}P_{\kappa,\delta,\ast}^{\left(  q\right)  }\left(  s^{\prime
}\right)  \ ,
\end{align*}
provided $q>3$ and $\lambda$ is sufficiently large, namely%
\begin{equation}
2^{\left(  \frac{6}{q}-2\right)  \lambda}+2^{\left(  \frac{3}{q}%
+\alpha-4\right)  \lambda}<\frac{1}{2C},\label{namely}%
\end{equation}
where $C$ is a constant arising from the early stages of the proof in
\textbf{Case 2}. Thus for $t>\lambda$ we have%
\begin{align*}
&  \max_{\lambda\leq s\leq t}P_{\kappa,\delta}^{\left(  q\right)  }\left(
s\right)  =\max_{\lambda\leq s\leq t}\sup_{\left\Vert f\right\Vert
_{L^{\infty}}\leq1}\left\Vert T\mathsf{Q}_{\kappa,s,U}f\right\Vert
_{L^{q}\left(  B_{\delta}\left(  0,2^{s}\right)  \right)  }\\
&  \lesssim C_{\varepsilon,N}2^{\left(  2+\alpha\right)  \lambda}2^{\frac
{6}{q}\lambda}\max_{\lambda<s\leq t}2^{\varepsilon s}+C_{\lambda,\kappa
,\delta,q}+\frac{1}{2}\max_{\lambda\leq s^{\prime}\leq s+10-\lambda}%
P_{\kappa,\delta}^{\left(  q\right)  }\left(  s^{\prime}\right)  +\frac{1}%
{2}\max_{\lambda<s^{\prime}\leq s+10-\lambda}P_{\kappa,\delta,\ast}^{\left(
q\right)  }\left(  s\right) \\
&  \lesssim C_{\varepsilon,N}2^{\left(  2+\alpha\right)  \lambda}2^{\frac
{6}{q}\lambda}\max_{\lambda<s\leq t}2^{\varepsilon s}+C_{\lambda,\kappa
,\delta,q}+\frac{1}{2}\max_{\lambda\leq s^{\prime}\leq s+10-\lambda}%
P_{\kappa,\delta}^{\left(  q\right)  }\left(  s\right)  ,
\end{align*}
and absorption now yields the inequality,%
\[
P_{\kappa,\delta}^{\left(  q\right)  }\left(  t\right)  \leq\max
_{\lambda<s\leq t}P_{\kappa,\delta}^{\left(  q\right)  }\left(  s\right)
\lesssim C_{\lambda,\kappa,\delta,q}+C_{\varepsilon}2^{\left(  2+\alpha
\right)  \lambda}2^{\frac{6}{q}\lambda}2^{\varepsilon t},\ \ \ \ \ \text{for
all }t>\lambda.
\]

If we take $\alpha=2$ and write $C=2^{c-1}$, then (\ref{namely}) becomes,%
\[
2^{-\frac{2}{q}\left(  q-3\right)  \lambda}+2^{-\frac{2}{q}\left(  q-\frac
{3}{2}\right)  \lambda}<\frac{1}{2C}=2^{-c},
\]
which is satisfied if both
\begin{align*}
\left(  \frac{2}{q}\left(  q-3\right)  \right)  \lambda &  >c+1\text{ and
}\left(  \frac{2}{q}\left(  q-\frac{3}{2}\right)  \right)  \lambda
>c+1\text{,}\\
\text{and in particular if }\lambda &  >\frac{c+1}{2}\frac{q}{q-3}\text{ and
}q>3,
\end{align*}
This shows that the `disjoint' parameter $\nu$ may be taken as
\begin{equation}
\nu<2^{10-\lambda}=2^{10}2^{-\frac{c+1}{2}\frac{q}{q-3}},\label{nu}%
\end{equation}
and so depends only on how much larger $q$ is than $3$. In conclusion, we have%
\begin{align*}
P_{\kappa,\delta}^{\left(  q\right)  }\left(  s\right)   &  \leq
C_{\lambda,\kappa,\delta,q}+C_{\varepsilon,\nu}2^{4\lambda}2^{2\left(
\lambda-10\right)  }2^{\varepsilon s}\leq C_{\lambda,\kappa,\delta
,q}+C_{\varepsilon,\nu}2^{6\lambda}2^{\varepsilon s}\\
&  =C_{\lambda,\kappa,\delta,q}+C_{\varepsilon,\nu}2^{60}\nu^{-6}%
2^{\varepsilon s}=C_{\lambda,\kappa,\delta,q}+C_{\varepsilon,\nu
}2^{\varepsilon s},
\end{align*}
for all $s>\lambda$ and $0<\nu<2^{10}2^{-\frac{c+1}{2}\frac{q}{q-3}}$, which
completes the proof of Theorem \ref{Loc lin} since%
\[
s>\lambda>\frac{c+1}{2}\frac{q}{q-3}\equiv C_{0}\frac{q}{q-3}.
\]

\end{proof}

\end{document}